\documentclass[leqno]{amsart}

\usepackage[leqno]{amsmath}
\usepackage{amssymb,amsthm}
\usepackage{enumerate, latexsym, mathrsfs, mdwlist}

\swapnumbers
\theoremstyle{definition}
\newtheorem{nul}{}[section]
\newtheorem{dfn}[nul]{Definition}

\newtheorem{cnstr}[nul]{Construction}
\newtheorem{ntn}[nul]{Notation}
\newtheorem{exm}[nul]{Example}

\newtheorem*{dfn*}{Definition}
\newtheorem*{axm*}{Axiom}
\newtheorem*{ntn*}{Notation}
\newtheorem*{exm*}{Example}
\newtheorem*{exr*}{Exercise}
\newtheorem*{int*}{Intuition}
\newtheorem*{qst*}{Question}

\theoremstyle{plain}

\newtheorem{thm}[nul]{Theorem}
\newtheorem{prp}[nul]{Proposition}
\newtheorem{lem}[nul]{Lemma}

\newtheorem{cor}{Corollary}[nul]
\newtheorem*{thm*}{Theorem}
\newtheorem*{prp*}{Proposition}
\newtheorem*{cor*}{Corollary}
\newtheorem*{lem*}{Lemma}
\newtheorem*{cnj*}{Conjecture}

\numberwithin{equation}{nul}

\DeclareMathOperator{\colim}{colim}

\DeclareMathOperator{\Fun}{Fun}
\DeclareMathOperator{\hocolim}{hocolim}
\DeclareMathOperator{\id}{id}

\DeclareMathOperator{\Map}{Map}
\DeclareMathOperator{\Mor}{Mor}

\DeclareMathOperator{\pr}{pr}

\newcommand{\FF}{\mathbf{F}}

\newcommand{\KK}{\mathbf{K}}

\newcommand{\QQ}{\mathbf{Q}}
\newcommand{\RR}{\mathbf{R}}
\renewcommand{\SS}{\mathbf{S}}
\newcommand{\TT}{\mathbf{T}}

\newcommand{\XX}{\mathbf{X}}

\newcommand{\Cat}{\mathbf{Cat}}

\newcommand{\Exact}{\mathbf{Exact}}

\newcommand{\Kan}{\mathbf{Kan}}

\newcommand{\Mod}{\mathbf{Mod}}

\newcommand{\Perf}{\mathbf{Perf}}

\newcommand{\Set}{\mathbf{Set}}

\newcommand{\VWald}{\mathbf{VWald}}
\newcommand{\Wald}{\mathbf{Wald}}

\newcommand{\cocart}{\mathrm{cocart}}

\newcommand{\op}{\mathrm{op}}

\newcommand{\coloneq}{\mathrel{\mathop:}=}

\makeatletter 
\def\revddots{\mathinner{\mkern1mu\raise\p@ 
\vbox{\kern7\p@\hbox{.}}\mkern2mu 
\raise4\p@\hbox{.}\mkern2mu\raise7\p@\hbox{.}\mkern1mu}} 
\makeatother

\usepackage{tikz}
\usetikzlibrary{matrix,arrows,decorations}

\pgfdeclaredecoration{single line}{initial}{\state{initial}[width=\pgfdecoratedpathlength-1sp]{\pgfpathmoveto{\pgfpointorigin}}\state{final}{\pgfpathlineto{\pgfpointorigin}}} 

\newcommand{\fromto}[2]{{#1}\ \tikz[baseline]\draw[>=stealth,->](0,0.5ex)--(0.5,0.5ex);\ {#2}}
\newcommand{\into}[2]{{#1}\ \tikz[baseline]\draw[>=stealth,right hook->](0,0.5ex)--(0.5,0.5ex);\ {#2}}

\newcommand{\cofto}[2]{{#1}\ \tikz[baseline]\draw[>=stealth,>->](0,0.5ex)--(0.5,0.5ex);\ {#2}}
\newcommand{\fibto}[2]{{#1}\ \tikz[baseline]\draw[>=stealth,->>](0,0.5ex)--(0.5,0.5ex);\ {#2}}
\newcommand{\equivto}[2]{{#1}\ \tikz[baseline]\draw[>=stealth,->,font=\scriptsize,inner sep=0.5pt](0,0.5ex)--node[above]{$\sim$}(0.5,0.5ex);\ {#2}}
\newcommand{\goesto}[2]{{#1}\ \tikz[baseline]\draw[|->](0,0.5ex)--(0.5,0.5ex);\ {#2}}

\renewcommand{\to}{\ \tikz[baseline]\draw[>=stealth,->](0,0.5ex)--(0.5,0.5ex);\ }
\newcommand{\ot}{\ \tikz[baseline]\draw[>=stealth,<-](0,0.5ex)--(0.5,0.5ex);\ }

\title{On the Q construction for exact $\infty$-categories}

\author{Clark Barwick}
\address{Massachusetts Institute of Technology, Department of Mathematics, Building 2, 77 Massachusetts Avenue, Cambridge, MA 02139-4307, USA}
\email{clarkbar@gmail.com}

\begin{document}

\begin{abstract} We prove that the algebraic $K$-theory of an exact $\infty$-category can be computed via an $\infty$-categorical variant of the $Q$ construction. This construction yields a quasicategory whose weak homotopy type is a delooping of the $K$-theory space. We show that the direct sum endows this homotopy type with the structure of a infinite loop space, which agrees with the canonical one. Finally, we prove a proto-d\'evissage result, which gives a necessary and sufficient condition for a \emph{nilimmersion} of stable $\infty$-categories to be a $K$-theory equivalence. In particular, we prove that a well-known conjecture of Ausoni--Rognes is equivalent to the weak contractibility of a particular $\infty$-category.
\end{abstract}

\maketitle

\tableofcontents

Exact $\infty$-categories, which we introduced in \cite{K21}, are a natural $\infty$-categorical generalization of Quillen's exact categories. They include a large portion of those $\infty$-categories to which one wishes to apply the machinery of Waldhausen's algebraic $K$-theory. The algebraic $K$-theory of any ordinary exact category, the $K$-theory of arbitrary schemes and stacks, and Waldhausen's $A$-theory of spaces can all be described as the $K$-theory of exact $\infty$-categories.

Quillen showed that the algberaic $K$-theory of an ordinary exact category can be described as the loopspace of the nerve of a category --- the $Q$ construction. In perfect analogy with this, we prove that the algebraic $K$-theory of an exact $\infty$-category can be computed as the loopspace of the classifying space of an $\infty$-category --- given by an $\infty$-categorical $Q$ construction. This construction yields a quasicategory whose weak homotopy type is a delooping of the $K$-theory space (Th. \ref{thm:Qworks}). Moreover, we show that the direct sum endows this homotopy type with the structure of a infinite loop space, which agrees with the canonical one (Th. \ref{prp:deloopwithoplus}). Finally, we discuss consequences of Quillen's Theorem A and Theorem B for $\infty$-categories (the latter of which we prove --- Th. \ref{thm:B}) for the algebraic $K$-theory of exact $\infty$-categories. In particular, we prove a ``proto-d\'evissage'' theorem (Th. \ref{thm:protodev}), which gives a necessary and sufficient condition for a \emph{nilimmersion} (Df. \ref{dfn:devissable}) of stable $\infty$-categories to be a $K$-theory equivalence. In particular, we show (Ex. \ref{exm:ARconj}) that the well-known conjecture of Ausoni--Rognes \cite[(0.2)]{MR1947457} can be expressed as the weak contractibility of a relatively simple $\infty$-category.


\section{Recollections on exact $\infty$-categories} We briefly recall the relevant definitions.

\begin{dfn} An $\infty$-category $C$ will be said to be \textbf{\emph{additive}} if its homotopy category $hC$ is additive (as a category enriched in the homotopy category $h\Kan$ of spaces).
\end{dfn}

\begin{exm} The nerve of any ordinary additive category is additive $\infty$-category, and any stable $\infty$-category is additive.
\end{exm}

\begin{dfn}\label{dfn:biWald} Suppose $\mathscr{C}$ an $\infty$-category, and suppose $\mathscr{C}_{\dag},\mathscr{C}^{\dag}\subset\mathscr{C}$ subcategories that contain all the equivalences. We call the morphisms of $\mathscr{C}_{\dag}$ \textbf{\emph{ingressive}}, and we call the morphisms of $\mathscr{C}^{\dag}$ \textbf{\emph{egressive}}.
\begin{enumerate}[(\ref{dfn:biWald}.1)]
\item A pullback square
\begin{equation*}
\begin{tikzpicture}[baseline]
\matrix(m)[matrix of math nodes, 
row sep=4ex, column sep=4ex, 
text height=1.5ex, text depth=0.25ex] 
{X&Y\\
X'&Y',\\}; 
\path[>=stealth,->,font=\scriptsize] 
(m-1-1) edge (m-1-2) 
edge (m-2-1) 
(m-1-2) edge[->>] (m-2-2) 
(m-2-1) edge[>->] (m-2-2); 
\end{tikzpicture}
\end{equation*}
is said to be \emph{ambigressive} if $\cofto{X'}{Y'}$ is ingressive and $\fibto{Y}{Y'}$ is egressive. Dually, a pushout square
\begin{equation*}
\begin{tikzpicture}[baseline]
\matrix(m)[matrix of math nodes, 
row sep=4ex, column sep=4ex, 
text height=1.5ex, text depth=0.25ex] 
{X&Y\\
X'&Y',\\}; 
\path[>=stealth,->,font=\scriptsize] 
(m-1-1) edge[>->] (m-1-2) 
edge[->>] (m-2-1) 
(m-1-2) edge (m-2-2) 
(m-2-1) edge (m-2-2); 
\end{tikzpicture}
\end{equation*}
is said to be \textbf{\emph{ambigressive}} if $\cofto{X}{Y}$ is ingressive and $\fibto{X}{X'}$ is egressive.
\item\label{item:exactinfty} We will say that the triple $(\mathscr{C},\mathscr{C}_{\dag},\mathscr{C}^{\dag})$ is an \textbf{\emph{exact $\infty$-category}} if it satisfies the following conditions.
\begin{enumerate}[(\ref{dfn:biWald}.\ref{item:exactinfty}.1)]
\item The underlying $\infty$-category $\mathscr{C}$ is additive.
\item The pair $(\mathscr{C},\mathscr{C}_{\dag})$ is a Waldhausen $\infty$-category.
\item The pair $(\mathscr{C},\mathscr{C}^{\dag})$ is a coWaldhausen $\infty$-category.
\item A square in $\mathscr{C}$ is an ambigressive pullback if and only if it is an ambigressive pushout.
\end{enumerate}
\item An \textbf{\emph{exact sequence}} in $\mathscr{C}$ is an ambigressive pushout/pullback square
\begin{equation*}
\begin{tikzpicture} 
\matrix(m)[matrix of math nodes, 
row sep=4ex, column sep=4ex, 
text height=1.5ex, text depth=0.25ex] 
{X'&X\\ 
0&X''\\}; 
\path[>=stealth,->,font=\scriptsize] 
(m-1-1) edge[>->] (m-1-2) 
edge[->>] (m-2-1) 
(m-1-2) edge[>->] (m-2-2) 
(m-2-1) edge[->>] (m-2-2); 
\end{tikzpicture}
\end{equation*}
in $\mathscr{C}$. The cofibration $\cofto{X'}{X}$ will be called the \textbf{\emph{fiber}} of $\fibto{X}{X''}$, and the fibration $\fibto{X}{X''}$ will be called the \textbf{\emph{cofiber}} of $\cofto{X'}{X}$.
\end{enumerate}
\end{dfn}

\begin{exm}\label{exm:exactinftycats}\begin{enumerate}[(\ref{exm:exactinftycats}.1)]
\item The nerve $NC$ of an ordinary category $C$ can be endowed with a triple structure yielding an exact $\infty$-category if and only if $C$ is an ordinary exact category, in the sense of Quillen, wherein the admissible monomorphisms are exactly the cofibrations, and the admissible epimorphisms are exactly the fibrations. This is proved by appealing to the ``minimal'' axioms of Keller \cite[App. A]{MR1052551}.
\item Any stable $\infty$-category is an exact $\infty$-category in which all morphisms are both egressive and ingressive.
\suspend{enumerate}
\noindent Suppose $\mathscr{A}$ a stable $\infty$-category equipped with a t-structure, and suppose $a$ and $b$ integers.
\resume{enumerate}[{[(\ref{exm:exactinftycats}.1)]}]
\item The $\infty$-category $\mathscr{A}_{[a,+\infty)}\coloneq\mathscr{A}_{\geq a}$ admits an exact $\infty$-category structure, in which every morphism is ingressive, but a morphism $\fromto{Y}{Z}$ is egressive just in case the induced morphism $\fromto{\pi_aX}{\pi_aY}$ is an epimorphism of $\mathscr{A}^{\heartsuit}$.
\item Dually, the $\infty$-category $\mathscr{A}_{(-\infty,b]}\coloneq\mathscr{A}_{\leq b}$ admits an exact $\infty$-category structure, in which every morphism is egressive, but a morphism $\fromto{X}{Y}$ is ingressive just in case the induced morphism $\fromto{\pi_bX}{\pi_bY}$ is a monomorphism of $\mathscr{A}^{\heartsuit}$.
\item We may intersect these subcategories to obtain the full subcategory
\begin{equation*}
\mathscr{A}_{[a,b]}\coloneq\mathscr{A}_{\geq a}\cap\mathscr{A}_{\leq b}\subset\mathscr{A},
\end{equation*}
and we may intersect the subcategories of ingressive and egressive morphisms described to obtain the following exact $\infty$-category structure on $\mathscr{A}_{[a,b]}$. A morphism $\fromto{X}{Y}$ is ingressive just in case the induced morphism $\fromto{\pi_bX}{\pi_bY}$ is a monomorphism of the abelian category $\mathscr{A}^{\heartsuit}$. A morphism $\fromto{Y}{Z}$ is egressive just in case the induced morphism $\fromto{\pi_aX}{\pi_aY}$ is an epimorphism of $\mathscr{A}^{\heartsuit}$.
\suspend{enumerate}
\noindent More generally, suppose $\mathscr{A}$ any stable $\infty$-category,
\resume{enumerate}[{[(\ref{exm:exactinftycats}.1)]}]
\item Suppose $\mathscr{C}\subset\mathscr{A}$ \emph{any} full additive subcategory that is closed under extensions. Declare a morphism $\fromto{X}{Y}$ of $\mathscr{C}$ to be ingressive just in case its cofiber in $\mathscr{A}$ lies in $\mathscr{C}$. Dually, declare a morphism $\fromto{Y}{Z}$ of $\mathscr{C}$ to be egressive just in case its fiber in $\mathscr{A}$ lies in $\mathscr{C}$. Then $\mathscr{C}$ is exact with this triple structure.
\end{enumerate}
\end{exm}

\begin{dfn} Suppose $\mathscr{C}$ and $\mathscr{D}$ two exact $\infty$-categories. A functor $\fromto{\mathscr{C}}{\mathscr{D}}$ that preserves both cofibrations and fibrations will be said to be \textbf{\emph{exact}} if both the functor
\begin{equation*}
\fromto{(\mathscr{C},\mathscr{C}_{\dag})}{(\mathscr{D},\mathscr{D}_{\dag})}
\end{equation*}
of Waldhausen $\infty$-categories and the functor
\begin{equation*}
\fromto{(\mathscr{C},\mathscr{C}^{\dag})}{(\mathscr{D},\mathscr{D}^{\dag})}
\end{equation*}
of coWaldhausen $\infty$-categories are exact.
\end{dfn}
\noindent It turns out that \emph{either} of these conditions suffices. See \cite[Pr. 3.4]{K21}.

Exact $\infty$-categories and exact functors between them organize themselves into an $\infty$-category $\Exact_{\infty}$.


\section{The twisted arrow $\infty$-category} In this section, we construct, for any exact $\infty$-category $\mathscr{C}$, an $\infty$-category $Q(\mathscr{C})$ whose underlying simplicial set is a (single) delooping of $K(\mathscr{C})$. The construction proceeds along the same principles as the ones used by Quillen; namely, we construct an $\infty$-category $Q(\mathscr{C})$ whose objects are the objects of $\mathscr{C}$, whose morphisms from $X$ to $Y$ are ``spans''
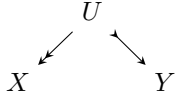
\begin{equation*}
\begin{tikzpicture} 
\matrix(m)[matrix of math nodes, 
row sep=3ex, column sep=3ex, 
text height=1.5ex, text depth=0.25ex] 
{&U&\\ 
X&&Y\\}; 
\path[>=stealth,->,font=\scriptsize] 
(m-1-2) edge[->>] (m-2-1) 
edge[>->] (m-2-3); 
\end{tikzpicture}
\end{equation*}
in which the morphism $\fibto{U}{X}$ is egressive and the morphism $\cofto{U}{Y}$ is ingressive, and whose composition law is given by the formation of pullbacks:
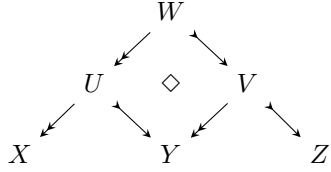
\begin{equation*}
\begin{tikzpicture} 
\matrix(m)[matrix of math nodes, 
row sep=3ex, column sep=3ex, 
text height=1.5ex, text depth=0.25ex] 
{&&W&&\\
&U&\Diamond&V&\\ 
X&&Y&&Z\\}; 
\path[>=stealth,->,font=\scriptsize] 
(m-1-3) edge[->>] (m-2-2) 
edge[>->] (m-2-4) 
(m-2-2) edge[->>] (m-3-1)
edge[>->] (m-3-3) 
(m-2-4) edge[->>] (m-3-3)
edge[>->] (m-3-5); 
\end{tikzpicture}
\end{equation*}
Some new machinery is involved in relating the $Q$ construction to our $\mathscr{S}$ construction. To codify their relationship, we form a Reedy fibrant straightening $\QQ_{\ast}(\mathscr{C})$ of an edgewise subdivision of the left fibration $\fromto{\iota_{N\Delta^{\op}}\mathscr{SC}}{N\Delta^{\op}}$; then we show that $\QQ_{\ast}(\mathscr{C})$ is in fact a \emph{complete Segal space} in the sense of Charles Rezk \cite{MR1804411}, whence by a theorem of Andr\'e Joyal and Myles Tierney \cite{JT}, the simplicial set $Q\mathscr{C}$, whose $n$-simplices may be identified with the vertices of the simplicial set $\QQ_n(\mathscr{C})$, is a quasicategory whose underlying simplicial set is equivalent to the geometric realization of $\QQ_{\ast}(\mathscr{C})$, which is in turn equivalent to $\iota_{N\Delta^{\op}}\mathscr{SC}$.

\begin{prp}\label{prp:subdivision} The following are equivalent for a functor $\theta\colon\fromto{\Delta}{\Delta}$.
\begin{enumerate}[(\ref{prp:subdivision}.1)]
\item The functor $\theta^{\op}\colon\fromto{N\Delta^{\op}}{N\Delta^{\op}}$ is cofinal in the sense of Joyal \cite[Df. 4.1.1.1]{HTT}.
\item The induced endofunctor $\theta^{\star}\colon\fromto{s\Set}{s\Set}$ on the ordinary category of simplicial sets (so that $(\theta^{\star}X)_n=X_{\theta(n)}$) carries every standard simplex $\Delta^m$ to a weakly contractible simplicial set.
\item The induced endofunctor $\theta^{\star}\colon\fromto{s\Set}{s\Set}$ on the ordinary category of simplicial sets is a left Quillen functor for the usual Quillen model structure.
\end{enumerate}
\begin{proof} By Joyal's variant of Quillen's Theorem A \cite[Th. 4.1.3.1]{HTT}, the functor $\theta^{\op}$ is cofinal just in case, for any integer $m\geq 0$, the nerve $N(\theta/\mathbf{m})$ is weakly contractible. The category $(\theta/\mathbf{m})$ is clearly equivalent to the category of simplices of $\theta^{\star}(\Delta^m)$, whose nerve is weakly equivalent to $\theta^{\star}(\Delta^m)$. This proves the equivalence of the first two conditions.

It is clear that for any functor $\theta\colon\fromto{\Delta}{\Delta}$, the induced functor $\theta^{\star}\colon\fromto{s\Set}{s\Set}$ preserves monomorphisms. Hence $\theta^{\star}$ is left Quillen just in case it preserves weak equivalences. Hence if $\theta^{\star}$ is left Quillen, then it carries the map $\equivto{\Delta^n}{\Delta^0}$ to an equivalence $\equivto{\theta^{\star}\Delta^n}{\theta^{\star}\Delta^0\cong\Delta^0}$, and, conversely, if $\theta^{\op}\colon\fromto{N\Delta^{\op}}{N\Delta^{\op}}$ is cofinal, then for any weak equivalence $\equivto{X}{Y}$, the induced map $\fromto{\theta^{\star}X}{\theta^{\star}Y}$ factors as
\begin{eqnarray}
\theta^{\star}X\simeq\hocolim_nX_{\theta(n)}&\simeq&\hocolim_nX_n\nonumber\\
&\simeq&X\nonumber\\
&\equivto{}{}&Y\nonumber\\
&\simeq&\hocolim_nY_n\nonumber\\
&\simeq&\hocolim_nY_{\theta(n)}\simeq\theta^{\star}Y,\nonumber
\end{eqnarray}
which is a weak equivalence. This proves the equivalence of the third condition with the first two.
\end{proof}
\end{prp}

\begin{dfn} Let us call any functor $\theta\colon\fromto{\Delta}{\Delta}$ satisfying the equivalent conditions above a \textbf{\emph{combinatorial subdivision}}.
\end{dfn}

\begin{nul} As pointed out to us by Katerina Velcheva, one may completely classify combinatorial subdivisions: they are generated by the functors $\id$ and $\op$ under the concatenation operation $\star$ on $\Delta$. Her work, which we hope will appear soon, also classifies more general sorts of subdivisions.
\end{nul}

Any combinatorial subdivision can be used to doctor our construction of $K$-theory in the following manner.
\begin{cnstr} Suppose $\theta\colon\fromto{\Delta}{\Delta}$ a combinatorial subdivision. Then pullback along $\theta$ defines an endofunctor $\theta^{\star}\colon\fromto{\Wald_{\infty/N\Delta^{\op}}^{\mathrm{cocart}}}{\Wald_{\infty/N\Delta^{\op}}^{\mathrm{cocart}}}$; since $\theta^{\op}$ is cofinal, it is clear that the following diagram commutes, up to equivalence:
\begin{equation*}
\begin{tikzpicture} 
\matrix(m)[matrix of math nodes, 
row sep=8ex, column sep=0ex, 
text height=1.5ex, text depth=0.25ex] 
{\Wald_{\infty/N\Delta^{\op}}^{\mathrm{cocart}}&&\Wald_{\infty/N\Delta^{\op}}^{\mathrm{cocart}}\\ 
&\VWald_{\infty}&\\}; 
\path[>=stealth,->,font=\scriptsize] 
(m-1-1) edge node[above]{$\theta^{\star}$} (m-1-3) 
edge node[below left]{$|\cdot|_{N\Delta^{\op}}$} (m-2-2) 
(m-1-3) edge node[below right]{$|\cdot|_{N\Delta^{\op}}$} (m-2-2); 
\end{tikzpicture}
\end{equation*}
In particular, the additivization of any pre-additive theory $\phi$ with left derived functor $\Phi$ can be computed as
\begin{equation*}
D\phi(\mathscr{C})\simeq\Omega\Phi|\theta^{\star}\mathscr{S}(\mathscr{C})|.
\end{equation*}
\end{cnstr}

The motivating example of a combinatorial subdivision is the following.

\begin{exm}\label{exm:edgewise} Denote by $\epsilon\colon\fromto{\Delta}{\Delta}$ the combinatorial subdivision given by the concatenation $\op\star\id$:
\begin{equation*}
\epsilon\colon\goesto{\mathbf{[n]}}{\mathbf{[n]}^{\op}\star\mathbf{[n]}\cong\mathbf{[2n+1]}}.
\end{equation*}
Including $\mathbf{[n]}$ into either factor of the join $\mathbf{[n]}^{\op}\star\mathbf{[n]}$ (either contravariantly or covariantly) defines two natural transformations $\fromto{\op}{\epsilon}$ and $\fromto{\id}{\epsilon}$. This functor induces an endofunctor $\epsilon^{\star}$ on the ordinary category of simplicial sets, together with natural transformations $\fromto{\epsilon^{\star}}{\op}$ and $\fromto{\epsilon^{\star}}{\id}$.

For any simplicial set $X$, the \textbf{\emph{edgewise subdivision}} of $X$ is the simplicial set
\begin{equation*}
\widetilde{\mathscr{O}}(X)\coloneq\epsilon^{\star}X.
\end{equation*}
That is, $\widetilde{\mathscr{O}}(X)$ is given by the formula
\begin{equation*}
\widetilde{\mathscr{O}}(X)_n=\Mor(\Delta^{n,\op}\star\Delta^n,X)\cong X_{2n+1}.
\end{equation*}
The two natural transformations described above give rise to a morphism
\begin{equation*}
\fromto{\widetilde{\mathscr{O}}(X)}{X^{\op}\times X},
\end{equation*}
functorial in $X$.
\end{exm}

\begin{nul}\label{nul:twarrofnerveisverve} For any simplicial set $X$, the vertices of $\widetilde{\mathscr{O}}(X)$ are edges of $X$; an edge of $\widetilde{\mathscr{O}}(X)$ from $\fromto{u}{v}$ to $\fromto{x}{y}$ can be viewed as a commutative diagram (up to chosen homotopy)
\begin{equation*}
\begin{tikzpicture} 
\matrix(m)[matrix of math nodes, 
row sep=6ex, column sep=6ex, 
text height=1.5ex, text depth=0.25ex] 
{u&x\\ 
v&y\\}; 
\path[>=stealth,->,font=\scriptsize] 
(m-1-2) edge (m-1-1) 
edge (m-2-2)
(m-1-1) edge (m-2-1)
(m-2-1) edge (m-2-2); 
\end{tikzpicture}
\end{equation*}
When $X$ is the nerve of an ordinary category $C$, $\widetilde{\mathscr{O}}(X)$ is isomorphic to the nerve of the twisted arrow category of $C$ in the sense of \cite{MR705421}. When $X$ is an $\infty$-category, we will therefore call $\widetilde{\mathscr{O}}(X)$ the \textbf{\emph{twisted arrow $\infty$-category}} of $X$. This terminology is justified by the following.
\end{nul}

\begin{prp}[Lurie, \protect{\cite[Pr. 4.2.3]{DAGX}}]\label{prp:twarrisinfincat} If $X$ is an $\infty$-category, then the functor $\fromto{\widetilde{\mathscr{O}}(X)}{X^{\op}\times X}$ is a left fibration; in particular, $\widetilde{\mathscr{O}}(X)$ is an $\infty$-category.
\end{prp}

\begin{exm} To illustrate, for any object $\mathbf{p}\in\Delta$, the $\infty$-category $\widetilde{\mathscr{O}}(\Delta^p)$ is the nerve of the category
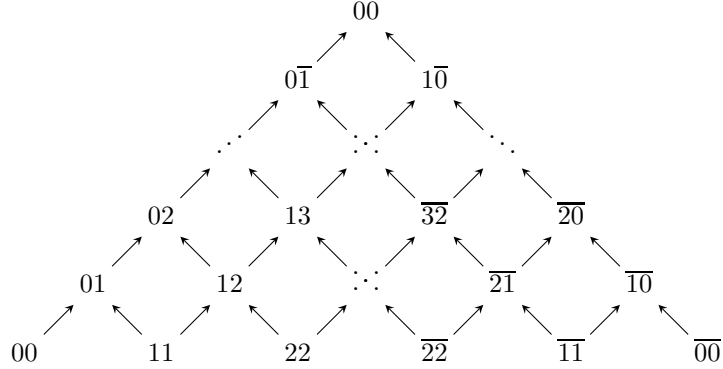
\begin{equation*}
\begin{tikzpicture} 
\matrix(m)[matrix of math nodes, 
row sep={6ex,between origins}, column sep={6ex,between origins}, 
text height=1.5ex, text depth=0.25ex] 
{&&&&&0\overline{0}&&&&&\\
&&&&0\overline{1}&&1\overline{0}&&&&\\
&&&\revddots&&\node{\ddots};\node{\revddots};&&\ddots&&&\\
&&02&&13&&\overline{3}\overline{2}&&\overline{2}\overline{0}&&\\
&01&&12&&\node{\ddots};\node{\revddots};&&\overline{2}\overline{1}&&\overline{1}\overline{0}&\\
00&&11&&22&&\overline{2}\overline{2}&&\overline{1}\overline{1}&&\overline{0}\overline{0}\\}; 
\path[>=stealth,->,font=\scriptsize] 
(m-2-5) edge (m-1-6) 
(m-2-7) edge (m-1-6)
(m-3-4) edge (m-2-5)
(m-3-6) edge (m-2-5)
(m-3-6) edge (m-2-7)
(m-3-8) edge (m-2-7)
(m-4-3) edge (m-3-4)
(m-4-5) edge (m-3-4)
(m-4-5) edge (m-3-6)
(m-4-7) edge (m-3-6)
(m-4-7) edge (m-3-8)
(m-4-9) edge (m-3-8)
(m-5-2) edge (m-4-3)
(m-5-4) edge (m-4-3)
(m-5-4) edge (m-4-5)
(m-5-6) edge (m-4-5)
(m-5-6) edge (m-4-7)
(m-5-8) edge (m-4-7)
(m-5-8) edge (m-4-9)
(m-5-10) edge (m-4-9)
(m-6-1) edge (m-5-2)
(m-6-3) edge (m-5-2)
(m-6-3) edge (m-5-4)
(m-6-5) edge (m-5-4)
(m-6-5) edge (m-5-6)
(m-6-7) edge (m-5-6)
(m-6-7) edge (m-5-8)
(m-6-9) edge (m-5-8)
(m-6-9) edge (m-5-10) 
(m-6-11) edge (m-5-10);
\end{tikzpicture}
\end{equation*}
(Here we write $\overline{n}$ for $p-n$.)
\end{exm}


\section{The $\infty$-categorical $Q$ construction} We now use the edgewise subdivision to define a quasicategorical variant of Quillen's $Q$ construction. To do this, it is convenient to use the theory of marked simplicial sets and the cocartesian model structure. As a result, in the definition and proposition that follow, we will use some of the notation of \cite{HTT}. 

\begin{dfn}\label{dfn:Rstar} For any marked simplicial set $X$, denote by $\RR_{\ast}(X)\colon\fromto{\Delta^{\op}}{s\Set}$ the functor given by the assignment
\begin{equation*}
\goesto{\mathbf{[n]}}{\Map^{\sharp}(\widetilde{\mathscr{O}}(\Delta^n)^{\op,\flat},X)}.
\end{equation*}
\end{dfn}

\begin{prp}\label{prp:Rstarisgood} The functor $\RR_{\ast}\colon\fromto{s\Set^{+}_{\cocart}}{\Fun(\Delta^{\op},s\Set)_{\mathrm{Reedy}}}$ preserves fibrant objects and weak equivalences between them.
\begin{proof} We first show that for any $\infty$-category $C$, the simplicial space $\RR_{\ast}(C^{\natural})$ is Reedy fibrant. This is the condition that for any monomorphism $\into{K}{L}$, the map
\begin{equation*}
\fromto{\Map^{\sharp}(\widetilde{\mathscr{O}}(L)^{\op,\flat},C^{\natural})}{\Map^{\sharp}(\widetilde{\mathscr{O}}(K)^{\op,\flat},C^{\natural})}
\end{equation*}
is a Kan fibration of simplicial sets. This follows immediately from Pr. \ref{prp:subdivision} and \cite[Lm. 3.1.3.6]{HTT}. To see that $\RR_{\ast}$ preserves weak equivalences between fibrant objects, we note that for any $\infty$-category $C$, the simplicial set $\Map^{\sharp}(\widetilde{\mathscr{O}}(\Delta^n)^{\op,\flat},C^{\natural})$ can be identified with the Kan complex $\iota\Fun(\widetilde{\mathscr{O}}(\Delta^n)^{\op},C)$, which clearly respects weak equivalences in $C$.
\end{proof}
\end{prp}

\begin{dfn}\label{dfn:ambigressfunct} Suppose $\mathscr{C}$ an exact $\infty$-category. For any integer $n\geq 0$, let us say that a functor $X\colon\fromto{\widetilde{\mathscr{O}}(\Delta^n)}{\mathscr{C}}$ is \textbf{\emph{ambigressive}} if, for any integers $0\leq i\leq k\leq\ell\leq j\leq n$, the square
\begin{equation*}
\begin{tikzpicture} 
\matrix(m)[matrix of math nodes, 
row sep=4ex, column sep=4ex, 
text height=1.5ex, text depth=0.25ex] 
{X_{ij}&X_{kj}\\ 
X_{i\ell}&X_{k\ell}\\}; 
\path[>=stealth,->,font=\scriptsize] 
(m-1-1) edge[>->] (m-1-2) 
edge[->>] (m-2-1) 
(m-1-2) edge[->>] (m-2-2) 
(m-2-1) edge[>->] (m-2-2); 
\end{tikzpicture}
\end{equation*}
is an ambigressive pullback. 

Write $\QQ_{\ast}(\mathscr{C})\subset\RR_{\ast}(\mathscr{C})$ for the subfunctor in which $\QQ_n(\mathscr{C})$ is the full simplicial subset of $\RR_n(\mathscr{C})$ spanned by the ambigressive functors $X\colon\fromto{\widetilde{\mathscr{O}}(\Delta^n)}{\mathscr{C}}$. Note that since any functor that is equivalent to an ambigressive functor is itself ambigressive, the simplicial set $\QQ_n(\mathscr{C})$ is a union of connected components of $\RR_n(\mathscr{C})$.
\end{dfn}

\begin{prp}\label{prp:QQstarisCSS} For any exact $\infty$-category $\mathscr{C}$, the simplicial space $\QQ_{\ast}(\mathscr{C})$ is a complete Segal space.
\begin{proof} The Reedy fibrancy of $\QQ_{\ast}(\mathscr{C})$ follows easily from the Reedy fibrancy of $\RR_{\ast}(\mathscr{C})$.

To see that $\QQ_{\ast}(\mathscr{C})$ is a Segal space, it is necessary to show that for any integer $n\geq 1$, the Segal map
\begin{equation*}
\fromto{\QQ_n(\mathscr{C})}{\QQ_1(\mathscr{C})\times_{\QQ_0(\mathscr{C})}\cdots\times_{\QQ_0(\mathscr{C})}\QQ_1(\mathscr{C})}
\end{equation*}
is an equivalence. Let $L_n$ denote the ordinary category
\begin{equation*}
00\ot 01\to 11\ot 12\to\cdots (n-1)(n-1)\ot (n-1)n\to nn;
\end{equation*}
equip $NL_n$ with the triple structure in which the maps $\fromto{(i-1)i}{ii}$ are ingressive, and the maps $\fromto{i(i+1)}{ii}$ are egressive. The target of the Segal map can then be identified with the maximal Kan complex contained in the full subcategory of $\Fun(NL_n,\mathscr{C})$ spanned by those functors $\fromto{NL_n}{\mathscr{C}}$ that preserve both cofibrations and fibrations. The Segal map is therefore an equivalence by the uniqueness of limits in $\infty$-categories \cite[Pr. 1.2.12.9]{HTT}.

Finally, to check that $\QQ_{\ast}(\mathscr{C})$ is complete, let $E$ be the nerve of the contractible ordinary groupoid with two objects; then completeness is equivalent to the assertion that the Rezk map
\begin{equation*}
\fromto{\QQ_0(\mathscr{C})}{\lim_{\mathbf{[n]}\in(\Delta/E)^{\op}}\QQ_n(\mathscr{C})}
\end{equation*}
is a weak equivalence. The source of this map can be identified with $\iota\mathscr{C}$; its target can be identified with the full simplicial subset of $\iota\Fun(\widetilde{\mathscr{O}}(E)^{\op},\mathscr{C})$ spanned by those functors $X\colon\fromto{\widetilde{\mathscr{O}}(E)^{\op}}{\mathscr{C}}$ such that for any simplex $\fromto{\Delta^n}{E}$, the induced functor $\fromto{\widetilde{\mathscr{O}}(\Delta^n)^{\op}}{\mathscr{C}}$ is ambigressive. Note that the twisted arrow category of the contractible ordinary groupoid with two objects is the contractible ordinary groupoid with four objects. Hence the image of any functor $X\colon\fromto{\widetilde{\mathscr{O}}(E)^{\op}}{\mathscr{C}}$ is contained in $\iota\mathscr{C}$, whence $X$ is automatically ambigressive. Thus the target of the Rezk map can be identified with $\iota\Fun(\widetilde{\mathscr{O}}(E)^{\op},\mathscr{C})$ itself, and the Rezk map is an equvalence.
\end{proof}
\end{prp}

Note that this result does not require the full strength of the condition that $\mathscr{C}$ be an exact $\infty$-category; it requires only that $\mathscr{C}$ admit a triple structure in which ambigressive pullbacks exist, the ambigressive pullback of an ingressive morphism is ingressive, and the ambigressive pullback of an egressive morphism is egressive.

\begin{nul} It is now clear that $\QQ_{\ast}$ defines a relative functor $\fromto{\Exact_{\infty}^0}{\mathbf{CSS}^0}$, where $\mathbf{CSS}^{\Delta}\subset\Fun(\Delta^{\op},s\Set)$ is the full simplicial subcategory spanned by complete Segal spaces (and $\mathbf{CSS}^0$ is its category of $0$-simplices). It therefore defines a functor of $\infty$-categories $\QQ_{\ast}\colon\fromto{\Exact_{\infty}}{\mathbf{CSS}}$. We may also regard $\QQ_{\ast}$ as a functor $\fromto{\Delta^{\op}\times\Exact_{\infty}^0}{s\Set}$.
\end{nul}

Now we aim to show that for an exact $\infty$-category $\mathscr{C}$, the simplicial space $\QQ_{\ast}(\mathscr{C})$ is a straightening of the left fibration $\fromto{\iota_{N\Delta^{\op}}\epsilon^{\star}\mathscr{SC}}{N\Delta^{\op}}$.

\begin{dfn}\label{dfn:tau} For any integer $n\geq 0$, denote by $\tau_n\colon\fromto{\widetilde{\mathscr{O}}(\Delta^n)^{\op}}{\mathscr{O}(\Delta^{n,\op}\star\Delta^n)}$ the fully faithful functor obtained as the nerve of the ordinary functor that carries an object $\fromto{i}{j}$ of the twisted arrow category of $\mathbf{[n]}$ to the object $\fromto{j}{i}$ of the arrow category of $\mathbf{[n]}^{\op}\star\mathbf{[n]}$, where $j$ is regarded as an object of $\mathbf{[n]}^{\op}$, and $i$ is regarded as an object of $\mathbf{[n]}$.

It is easy to verify that the functors $\tau_n$ are compatible with all face and degeneracy maps, so they fit together to form a natural transformation $\tau\colon\fromto{\op\circ\widetilde{\mathscr{O}}}{\mathscr{O}\circ\epsilon}$ of functors $\fromto{\Delta}{\Cat_{\infty}}$.
\end{dfn}

Once again we employ the uniqueness of limits and colimits in $\infty$-categories \cite[Pr. 1.2.12.9]{HTT} to deduce the following.
\begin{prp}\label{prp:taugivesequiv} The natural transformation $\tau$ induces an equivalence
\begin{equation*}
\equivto{\iota\circ\epsilon^{\star}\widetilde{\SS}_{\ast}}{\TT_{\ast}},
\end{equation*}
where $\widetilde{\SS}_{\ast}$ is the functor described in \cite[4.16]{K21}.
\end{prp}

\begin{cor} There is a natural equivalence $\equivto{K(\mathscr{C})}{\Omega|\QQ_{\ast}(\mathscr{C})|}$ for any exact $\infty$-category $\mathscr{C}$.
\end{cor}

Joyal and Tierney show that the functor that carries a simplicial space $\XX$ to the simplicial set whose $n$-simplices are the vertices of $\XX_n$ induces an equivalence of relative categories $\fromto{\mathbf{CSS}^0}{\Cat_{\infty}^0}$. This leads us to the following definition and theorem.

\begin{dfn} For any exact $\infty$-category $\mathscr{C}$, denote by $Q(\mathscr{C})$ the $\infty$-category whose $n$-simplices are vertices of $\QQ_n(\mathscr{C})$, i.e., ambigressive functors $\fromto{\widetilde{\mathscr{O}}(\Delta^n)^{\op}}{\mathscr{C}}$. This defines a relative functor $Q\colon\fromto{\Exact_{\infty}^0}{\Cat_{\infty}^0}$ and hence a functor of $\infty$-categories $Q\colon\fromto{\Exact_{\infty}}{\Cat_{\infty}}$.
\end{dfn}

\begin{nul} For any exact $\infty$-category, an $n$-simplex of $Q(\mathscr{C})$ is a diagram
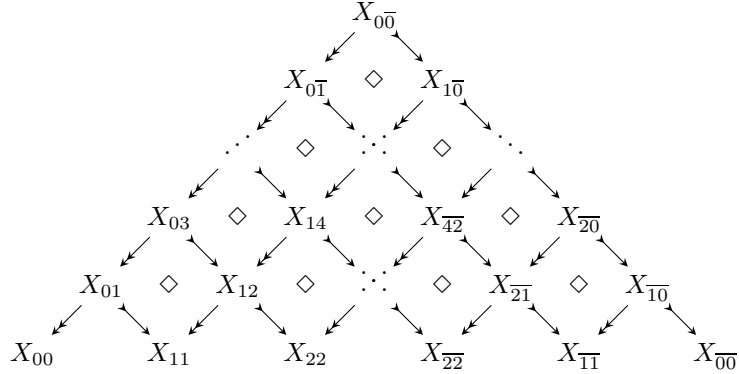
\begin{equation*}
\begin{tikzpicture} 
\matrix(m)[matrix of math nodes, 
row sep={6ex,between origins}, column sep={6ex,between origins}, 
text height=1.5ex, text depth=0.25ex] 
{&&&&&X_{0\overline{0}}&&&&&\\
&&&&X_{0\overline{1}}&\Diamond&X_{1\overline{0}}&&&&\\
&&&\revddots&\Diamond&\node{\ddots};\node{\revddots};&\Diamond&\ddots&&&\\
&&X_{03}&\Diamond&X_{14}&\Diamond&X_{\overline{4}\overline{2}}&\Diamond&X_{\overline{2}\overline{0}}&&\\
&X_{01}&\Diamond&X_{12}&\Diamond&\node{\ddots};\node{\revddots};&\Diamond&X_{\overline{2}\overline{1}}&\Diamond&X_{\overline{1}\overline{0}}&\\
X_{00}&&X_{11}&&X_{22}&&X_{\overline{2}\overline{2}}&&X_{\overline{1}\overline{1}}&&X_{\overline{0}\overline{0}}\\}; 
\path[>=stealth,<-,font=\scriptsize] 
(m-2-5) edge[<<-] (m-1-6) 
(m-2-7) edge[<-<] (m-1-6)
(m-3-4) edge[<<-] (m-2-5)
(m-3-6) edge[<-<] (m-2-5)
(m-3-6) edge[<<-] (m-2-7)
(m-3-8) edge[<-<] (m-2-7)
(m-4-3) edge[<<-] (m-3-4)
(m-4-5) edge[<-<] (m-3-4)
(m-4-5) edge[<<-] (m-3-6)
(m-4-7) edge[<-<] (m-3-6)
(m-4-7) edge[<<-] (m-3-8)
(m-4-9) edge[<-<] (m-3-8)
(m-5-2) edge[<<-] (m-4-3)
(m-5-4) edge[<-<] (m-4-3)
(m-5-4) edge[<<-] (m-4-5)
(m-5-6) edge[<-<] (m-4-5)
(m-5-6) edge[<<-] (m-4-7)
(m-5-8) edge[<-<] (m-4-7)
(m-5-8) edge[<<-] (m-4-9)
(m-5-10) edge[<-<] (m-4-9)
(m-6-1) edge[<<-] (m-5-2)
(m-6-3) edge[<-<] (m-5-2)
(m-6-3) edge[<<-] (m-5-4)
(m-6-5) edge[<-<] (m-5-4)
(m-6-5) edge[<<-] (m-5-6)
(m-6-7) edge[<-<] (m-5-6)
(m-6-7) edge[<<-] (m-5-8)
(m-6-9) edge[<-<] (m-5-8)
(m-6-9) edge[<<-] (m-5-10) 
(m-6-11) edge[<-<] (m-5-10);
\end{tikzpicture}
\end{equation*}
of $\mathscr{C}$ in which every square is an ambigressive pullback/pushout. (Here we write $\overline{n}$ for $p-n$.) 
\end{nul}

\begin{thm}\label{thm:Qworks} There is a natural equivalence $\equivto{K(\mathscr{C})}{\Omega Q(\mathscr{C})}$ for any exact $\infty$-category $\mathscr{C}$.
\end{thm}

As a final note, let us note that the $\infty$-categorical $Q$ construction we have introduced here is an honest generalization of Quillen's original $Q$ construction.
\begin{prp} If $C$ is an ordinary exact category, then $Q(NC)$ is canonically equivalent to the nerve $N(QC)$ of Quillen's $Q$ construction.
\begin{proof} Unwinding the definitions, we find that the simplicial set $\mathrm{Hom}^R_{Q(NC)}(X,Y)$ of \cite[\S 1.2.2]{HTT} is isomorphic to the nerve of the groupoid $G(X,Y)$ in which an object $Z$ is a diagram
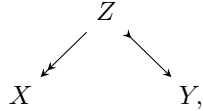
\begin{equation*}
\begin{tikzpicture} 
\matrix(m)[matrix of math nodes, 
row sep=4ex, column sep=4ex, 
text height=1.5ex, text depth=0.25ex] 
{&Z&\\ 
X&&Y,\\}; 
\path[>=stealth,->,font=\scriptsize] 
(m-1-2) edge[->>] (m-2-1) 
edge[>->] (m-2-3); 
\end{tikzpicture}
\end{equation*}
where $\fibto{Z}{X}$ is an admissible epimorphism and $\cofto{Z}{Y}$ is an admissible monomorphism, and in which a morphism $\fromto{Z'}{Z}$ is a diagram
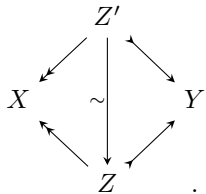
\begin{equation*}
\begin{tikzpicture} 
\matrix(m)[matrix of math nodes, 
row sep=4ex, column sep=4ex, 
text height=1.5ex, text depth=0.25ex] 
{&Z'&\\
X&&Y\\
&Z&.\\}; 
\path[>=stealth,->,font=\scriptsize] 
(m-1-2) edge[->>] (m-2-1)
edge[inner sep=0.5pt] node[left]{$\sim$} (m-3-2)
edge[>->] (m-2-3)
(m-3-2) edge[->>] (m-2-1) 
edge[>->] (m-2-3); 
\end{tikzpicture}
\end{equation*}
It is now immediate that $N(QC)$ is equivalent to the homotopy category of $Q(NC)$, and so it suffices by \cite[Pr. 2.3.4.18]{HTT} to verify that every connected component of the groupoid $G(X,Y)$ is contractible. For this, we simply note that if $f,g$ are two morphisms $\fromto{Z'}{Z}$ of $G(X,Y)$, then $f=g$ since $\cofto{Z}{Y}$ is a monomorphism.
\end{proof}
\end{prp}


\section{An infinite delooping of $Q$} In this section, we describe an infinite delooping of the $Q$ construction of the previous section. In effect, the direct sum on an exact $\infty$-category $\mathscr{C}$ induces a symmetric monoidal structure on the $\infty$-category $Q(\mathscr{C})$. Segal's delooping machine then applies to give an infinite delooping of the $Q$ construction. We will then show that this delooping coincides with the canonical one from \cite[Cor. 7.4.1]{K1} by appealing to the Additivity Theorem \cite[Th. 7.2]{K1}. We thank Lars Hesselholt for his questions about such a delooping; his hunch was right all along.

It is convenient to introduce some ordinary categories that control the combinatorics of direct sums.

\begin{ntn} Denote by $\Lambda(\FF)$ the following ordinary category. An object of $\Lambda(\FF)$ is a finite set; a morphism $\fromto{J}{I}$ of $\Lambda(\FF)$ is a map $\fromto{J}{I_{+}}$, or equivalently a pointed map $\fromto{J_{+}}{I_{+}}$. Clearly $\Lambda(\FF)$ is isomorphic to the category of pointed finite sets, but we shall regard the objects of $\Lambda(\FF)$ simply as finite (unpointed) sets.

Recall that $\Lambda(\FF)$ admits a symmetric monoidal structure that carries a pair of finite sets to their product. We shall denote this symmetric monoidal structure $\goesto{(I,J)}{I\wedge J}$.

For any finite set $I$, denote by $L_I$ the following ordinary category. An object of $L_I$ is a subset $J\subset I$. A morphism from an object $K\subset I$ to an object $J\subset I$ is a map $\psi\colon\fromto{K}{J_+}$ such that the square
\begin{equation*}
\begin{tikzpicture} 
\matrix(m)[matrix of math nodes, 
row sep=4ex, column sep=4ex, 
text height=1.5ex, text depth=0.25ex] 
{\psi^{-1}(J)&J\\ 
K&I\\}; 
\path[>=stealth,->,font=\scriptsize] 
(m-1-1) edge (m-1-2) 
edge (m-2-1) 
(m-1-2) edge[right hook->] (m-2-2) 
(m-2-1) edge[right hook->] (m-2-2); 
\end{tikzpicture}
\end{equation*}
commutes.
\begin{equation*}
\end{equation*}

Any morphism $\phi\colon\fromto{I^{\,\prime}}{I}$ of $\Lambda(\FF)$ induces a functor $\phi^{\star}\colon\fromto{L_{I^{\,\prime}}}{L_I}$ that carries $J\subset I$ to $\phi^{-1}(J)\subset I^{\,\prime}$. With this, one confirms easily that the assignments $\goesto{I}{L_I}$ and $\goesto{\phi}{\phi^{\star}}$ together define a functor $L\colon\fromto{\Gamma}{\Cat}$.

If $K\subset J$, then write $i_{K\subset J}$ for the morphism $\fromto{K}{J}$ of $\Lambda(\FF)$ given by the inclusion $\into{K}{J_{+}}$, and write $p_{K\subset J}$ for the morphism $\fromto{J}{K}$ of $\Lambda(\FF)$ given by the map $\into{J}{K_{+}}$ that carries every $j\in J\setminus K$ to the basepoint, and every element $j\in K$ to itself.
\end{ntn}

We now proceed to show that an exact $\infty$-category $\mathscr{C}$ admits an essentially unique symmetric monoidal structure in which all the multiplication functors are exact.

\begin{ntn} We write $\Exact_{\infty}^{\oplus}$ for the full subcategory of $\Fun(N\Lambda(\FF),\Exact_{\infty})$ consisting of those functors $\mathscr{X}\colon\fromto{N\Lambda(\FF)}{\Exact_{\infty}}$ such that for any finite set $I$, the functors $\{\fromto{\mathscr{X}(I)}{\mathscr{X}(\{i\})}\}_{i\in I}$ induced by the morphisms $p_{\{i\}\subset I}\colon\fromto{I}{\{i\}}$ together exhibit the exact $\infty$-category $\mathscr{X}(I)$ as a product of the exact $\infty$-categories $\mathscr{X}(\{i\})$. Write $U\colon\fromto{\Exact_{\infty}^{\oplus}}{\Exact_{\infty}}$ for the evaluation functor $\goesto{\mathscr{X}}{\mathscr{X}(\{1\})}$.

In the opposite direction, let us construct a functor $\mathbf{DS}\colon\fromto{\Exact_{\infty}}{\Exact_{\infty}^{\oplus}}$. For any finite set $I$ and any exact $\infty$-category, denote by $\mathbf{DS}(I;\mathscr{C})$ the full subcategory of $\Fun(NL_I,\mathscr{C})$ spanned by those functors $X\colon\fromto{NL_I}{\mathscr{C}}$ such that, for any subset $J\subset I$,
\begin{enumerate}[(\ref{prp:exactdirsum}.1)]
\item the set of morphisms
\begin{equation*}
\{\fromto{X(J)}{X(\{j\})}\}_{j\in J}
\end{equation*}
induced by the morphisms $p_{\{j\}\subset J}\colon\fromto{J}{\{j\}}$ of $L_I$ exhibit $X(J)$ as a product of the objects $X(\{j\})$, and,
\item dually, the morphisms $\{\fromto{X(\{j\})}{X(J)}\}_{j\in J}$ induced by the morphisms $i_{\{j\}\subset J}\colon\fromto{\{j\}}{J}$ of $L_I$ exhibit $X(J)$ as a coproduct of the objects $X(\{j\})$.
\end{enumerate}
It is clear that $\goesto{(I,\mathscr{C})}{\mathbf{DS}(I;\mathscr{C})}$ defines a simplicial functor
\begin{equation*}
\fromto{N\Lambda(\FF)\times\Exact_{\infty}^{\Delta}}{\Cat_{\infty}^{\Delta}}.
\end{equation*}

Now suppose $\mathscr{C}$ an exact $\infty$-category. It follows from the uniqueness of limits and colimits in $\infty$-categories \cite[Pr. 1.2.12.9]{HTT} that for any finite set $I$, the functors
\begin{equation*}
\{\fromto{\mathbf{DS}(I;\mathscr{C})}{\mathbf{DS}(\{i\},\mathscr{C})}\}_{i\in I}
\end{equation*}
induced by the morphisms $p_{\{i\}\subset I}\colon\fromto{I}{\{i\}}$ of $\Lambda(\FF)$ together exhibit the $\infty$-category $\mathbf{DS}(I;\mathscr{C})$ as the product of the $\infty$-categories $\mathbf{DS}(\{i\},\mathscr{C})$, which are each in turn equivalent to $\mathscr{C}$. Consequently the $\infty$-categories $\mathbf{DS}(I;\mathscr{C})$ are exact $\infty$-categories, and since direct sum in $\mathscr{C}$ preserves ingressives and any pushouts that exist, one easily confirms that $\mathbf{DS}(-;\mathscr{C})$ is an object of $\Exact_{\infty}^{\oplus}$; this therefore defines a functor $\mathbf{DS}\colon\fromto{\Exact_{\infty}}{\Exact_{\infty}^{\oplus}}$.
\end{ntn}

\begin{prp}\label{prp:exactdirsum} The functor $U$ exhibits an equivalence $\equivto{\Exact_{\infty}^{\oplus}}{\Exact_{\infty}}$, and the functor $\mathbf{DS}$ exhibits a quasi-inverse to it.
\begin{proof} In light of the discussion above, it suffices to prove that for any object $\mathscr{X}\in\Exact_{\infty}^{\oplus}$, the cocartesian fibration $\fromto{\mathscr{X}^{\oplus}}{N\Lambda(\FF)}$ classified by the composite
\begin{equation*}
N\Lambda(\FF)\ \tikz[baseline]\draw[>=stealth,->](0,0.5ex)--(0.5,0.5ex);\ \Exact_{\infty}\ \tikz[baseline]\draw[>=stealth,->](0,0.5ex)--(0.5,0.5ex);\ \Cat_{\infty}
\end{equation*}
is a cartesian symmetric monoidal $\infty$-category \cite[Df. 2.4.0.1]{HA}. It is obvious that it is symmetric monoidal. Furthermore, since the map $\fromto{\Delta^0\simeq\mathscr{X}(\varnothing)}{\mathscr{X}(\{1\})}$ is exact, the unit object is a zero object, and since the functor
\begin{equation*}
\otimes\colon\fromto{\mathscr{X}(\{1\})\times\mathscr{X}(\{1\})\simeq\mathscr{X}(\{1,2\})}{\mathscr{X}(\{1\})}
\end{equation*}
induced by the unique map $\fromto{\{1,2\}}{\{1\}_{+}}$ that does not hit the basepoint is exact, when applied to pullback squares
\begin{equation*}
\begin{tikzpicture}[baseline]
\matrix(m)[matrix of math nodes, 
row sep=4ex, column sep=4ex, 
text height=1.5ex, text depth=0.25ex] 
{X&0\\ 
X&0\\}; 
\path[>=stealth,->>,font=\scriptsize] 
(m-1-1) edge (m-1-2) 
edge (m-2-1) 
(m-1-2) edge (m-2-2) 
(m-2-1) edge (m-2-2); 
\end{tikzpicture}
\textrm{\quad and\quad}
\begin{tikzpicture}[baseline]
\matrix(m)[matrix of math nodes, 
row sep=4ex, column sep=4ex, 
text height=1.5ex, text depth=0.25ex] 
{Y&Y\\ 
0&0,\\}; 
\path[>=stealth,->>,font=\scriptsize] 
(m-1-1) edge (m-1-2) 
edge (m-2-1) 
(m-1-2) edge (m-2-2) 
(m-2-1) edge (m-2-2); 
\end{tikzpicture}
\end{equation*}
it yields a pullback square
\begin{equation*}
\begin{tikzpicture} 
\matrix(m)[matrix of math nodes, 
row sep=4ex, column sep=4ex, 
text height=1.5ex, text depth=0.25ex] 
{X\otimes Y&Y\\ 
X&0.\\}; 
\path[>=stealth,->>,font=\scriptsize] 
(m-1-1) edge (m-1-2) 
edge (m-2-1) 
(m-1-2) edge (m-2-2) 
(m-2-1) edge (m-2-2); 
\end{tikzpicture}
\end{equation*}
Thus $\mathscr{X}^{\oplus}$ is cartesian. 
\end{proof}
\end{prp}

\begin{nul} We may now lift the $Q$ construction to $\Exact_{\infty}^{\oplus}$. Composition with the functors
\begin{equation*}
Q\colon\fromto{\Exact_{\infty}}{\Cat_{\infty}}\textrm{\quad and\quad}\QQ_{\ast}\colon\fromto{\Exact_{\infty}}{\mathbf{CSS}}
\end{equation*}
defines functors
\begin{equation*}
\fromto{\Fun(N\Lambda(\FF),\Exact_{\infty})}{\Fun(N\Lambda(\FF),\Cat_{\infty})}
\end{equation*}
and
\begin{equation*}
\fromto{\Fun(N\Lambda(\FF),\Exact_{\infty})}{\Fun(N\Lambda(\FF),\mathbf{CSS})},
\end{equation*}
which restrict to functors
\begin{equation*}
Q^{\oplus}\colon\fromto{\Exact_{\infty}^{\oplus}}{\Cat_{\infty}^{\otimes}}\textrm{\quad and\quad}\QQ_{\ast}^{\oplus}\colon\fromto{\Exact_{\infty}^{\oplus}}{\mathbf{CSS}^{\otimes}}
\end{equation*}
where $\Cat_{\infty}^{\otimes}$ (respectively, $\mathbf{CSS}^{\otimes}$) denotes the full subcategory of the $\infty$-category $\Fun(N\Lambda(\FF),\Cat_{\infty})$ (respectively, of $\Fun(N\Lambda(\FF),\mathbf{CSS})$) spanned by those functors $X\colon\fromto{N\Lambda(\FF)}{\Cat_{\infty}}$ (respectively, by those functors $X\colon\fromto{N\Lambda(\FF)}{\mathbf{CSS}}$) such that for any finite set $I$, the morphisms $\{\fromto{X(I)}{X(\{i\})}\}_{i\in I}$ induced by the morphisms $p_{\{i\}\subset I}\colon\fromto{I}{\{i\}}$ together exhibit $X(I)$ as a product of $X(\{i\})$. We thus have commutative squares
\begin{equation*}
\begin{tikzpicture}[baseline]
\matrix(m)[matrix of math nodes, 
row sep=4ex, column sep=4ex, 
text height=1.5ex, text depth=0.25ex] 
{\Exact_{\infty}^{\oplus}&\Cat_{\infty}^{\otimes}\\ 
\Exact_{\infty}&\Cat_{\infty}\\}; 
\path[>=stealth,->,font=\scriptsize] 
(m-1-1) edge node[above]{$Q^{\oplus}$} (m-1-2) 
edge node[left]{$U$} (m-2-1) 
(m-1-2) edge node[right]{$U$} (m-2-2) 
(m-2-1) edge node[below]{$Q$}(m-2-2); 
\end{tikzpicture}
\textrm{\quad and\quad}
\begin{tikzpicture}[baseline]
\matrix(m)[matrix of math nodes, 
row sep=4ex, column sep=4ex, 
text height=1.5ex, text depth=0.25ex] 
{\Exact_{\infty}^{\oplus}&\mathbf{CSS}^{\otimes}\\ 
\Exact_{\infty}&\mathbf{CSS}\\}; 
\path[>=stealth,->,font=\scriptsize] 
(m-1-1) edge node[above]{$\QQ_{\ast}^{\oplus}$} (m-1-2) 
edge node[left]{$U$} (m-2-1) 
(m-1-2) edge node[right]{$U$} (m-2-2) 
(m-2-1) edge node[below]{$\QQ_{\ast}$}(m-2-2); 
\end{tikzpicture}
\end{equation*}
in which the vertical maps are evaluation at the vertex $\{1\}\in\Lambda(\FF)$.

For any object $\mathscr{X}\in\Exact_{\infty}^{\oplus}$, the object $Q^{\oplus}\mathscr{X}\simeq|\QQ_{\ast}^{\oplus}\mathscr{X}|$ is in particular a $\Gamma$-space that satisfies the Segal condition and, since the underlying space of $Q\mathscr{X}(\{1\})$ is connected, it is grouplike. Segal's delooping machine therefore provides an infinite delooping of $Q\mathscr{X}(\{1\})$.
\end{nul}

Let us now prove that the infinite delooping provided by $Q^{\oplus}$ agrees with the one guaranteed by \cite[Cor. 7.4.1]{K1}. Roughly speaking, the approach is the following. For any $\Gamma$-space $X$ that satisfies the Segal condition that is grouplike, the $n$-fold delooping of $X(\{1\})$ can be obtained as the geometric realization of $X$ when precomposed with a certain multisimplicial object $\fromto{(\Delta^{\op})^{\times n}}{\Lambda(\FF)}$. We now recall the construction of that multisimplicial object, and we show that when $X=Q^{\oplus}\mathscr{C}$ for an exact $\infty$-category $\mathscr{C}$, the $n$-fold delooping obtained in this way coincides with the iterated $\mathscr{S}$ construction.

\begin{ntn}\label{ntn:functoru} Write $u\colon\fromto{\Delta^{\op}}{\Lambda(\FF)}$ for the following functor. For any nonempty totally ordered finite set $S$, let $u(S)$ be the set of surjective morphisms $\fromto{S}{\mathbf{[1]}}$ of $\Delta$. For any map $g\colon\fromto{S}{T}$ of $\Delta$, define the map $u(g)\colon\fromto{u(T)}{u(S)_{+}}$ by
\begin{equation*}
u(g)(\eta)=\begin{cases}
\eta\circ g&\textrm{if $\eta\circ g$ is surjective;}\\
\ast&\textrm{otherwise.}
\end{cases}
\end{equation*}
Furthermore, for any nonnegative integer $n$, we may define $u^{(n)}\colon\fromto{(\Delta^{\op})^{\times n}}{\Lambda(\FF)}$ as the composite
\begin{equation*}
(\Delta^{\op})^{\times n}\ \tikz[baseline]\draw[>=stealth,->,font=\scriptsize](0,0.5ex)--node[above]{$u^n$}(0.5,0.5ex);\ (\Lambda(\FF))^{\times n}\ \tikz[baseline]\draw[>=stealth,->,font=\scriptsize](0,0.5ex)--node[above]{$\wedge^n$}(0.5,0.5ex);\ \Lambda(\FF).
\end{equation*}

For any nonempty totally ordered finite set $S$ and any morphism $\alpha\colon\fromto{\mathbf{[1]}}{S}$ of $\Delta$, let $\rho_S(\alpha)\subset u(S)$ be the set of retractions of $\alpha$, i.e., the morphisms $\beta\colon\fromto{S}{\mathbf{[1]}}$ such that $\beta\circ\alpha=\id$. For any nonnegative integer $n$, denote by $O\colon\fromto{\Delta^{\times n}}{\Cat}$ for the multicosimplicial category given by
\begin{equation*}
O(S_1,S_2,\dots,S_n)\coloneq\Fun(\mathbf{[1]},S_1)\times\Fun(\mathbf{[1]},S_2)\times\cdots\times\Fun(\mathbf{[1]},S_n),
\end{equation*}
and write $\rho\colon\fromto{O}{u^{(n),\op,\star}L}$ for the natural transformation given by
\begin{equation*}
\rho_{S_1,S_2,\dots,S_n}(\alpha_1,\alpha_2,\dots,\alpha_n)\coloneq\rho_{S_1}(\alpha_1)\wedge\rho_{S_2}(\alpha_2)\wedge\cdots\wedge\rho_{S_n}(\alpha_n).
\end{equation*}
This natural transformation induces a natural transformation
\begin{equation*}
\rho^{\star}\colon\fromto{u^{(n),\op}\mathbf{DS}}{\widetilde{\SS}^{(n)}_{\ast}}
\end{equation*}
between functors $\fromto{(\Delta^{\op})^{\times n}\times\Exact_{\infty}}{\Wald_{\infty}}$. Combined with the natural transformation $\tau$ from Df. \ref{dfn:tau}, we obtain a natural transformation
\begin{equation*}
R\colon\fromto{u^{(n),\star}\QQ_{\ast}^{\oplus}\circ\mathbf{DS}}{\iota\circ\epsilon^{\star}\widetilde{\SS}_{\ast}\circ\widetilde{\SS}_{\ast}^{(n)}}
\end{equation*}
between functors $\fromto{\Delta^{\op}\times(\Delta^{\op})^{\times n}\times\Exact_{\infty}}{\Wald_{\infty}}$
\end{ntn}

The following is an immediate consequence of the Additivity Theorem \cite[Th. 7.2]{K1}.
\begin{prp}\label{prp:deloopwithoplus} For any positive integer $n$, the natural transformation $R$ above is a natural equivalence
\begin{equation*}
|u^{(n),\star}\QQ_{\ast}^{\oplus}\circ\mathbf{DS}|\simeq|\iota\circ\epsilon^{\star}\widetilde{\SS}_{\ast}\circ\widetilde{\SS}_{\ast}^{(n)}|
\end{equation*}
of functors $\fromto{\Exact_{\infty}}{\Kan}$.
\end{prp}
\noindent In other words, the delooping of algebraic $K$-theory provided by the $Q^{\oplus}$ construction is naturally equivalent to the canonical delooping obtained in \cite[Cor. 7.4.1]{K1}.


\section{Nilimmersions and a relative $Q$ construction} One of the original uses of the $Q$ construction was Quillen's proof of the D\'evissage Theorem, which gives an effective way of determining whether an exact functor $\psi\colon\fromto{\mathscr{B}}{\mathscr{A}}$ induces a $K$-theory equivalence \cite[Th. 4]{MR0338129}. This theorem made it possible for Quillen to identify the ``fiber term'' in his Localization Sequence for higher algebraic $K$-theory \cite[Th. 5]{MR0338129}. The technical tool Quillen introduced for this purpose was his celebrated Theorem A. Joyal proved the following $\infty$-categorical variant of Quillen's Theorem A.
\begin{thm}[Theorem A for $\infty$-categories \protect{\cite[Th. 4.1.3.1]{HTT}}] Suppose $G\colon\fromto{C}{D}$ a functor between $\infty$-categories. If, for any object $X\in D$, the $\infty$-category
\begin{equation*}
G_{X/}\coloneq D_{X/}\times_DC
\end{equation*}
is weakly contractible, then the map of simplicial sets $G$ is a weak homotopy equivalence.
\end{thm}
\noindent This form of Theorem A (or rather its opposite) directly implies a recognition principle for $K$-theory equivalences in the $\infty$-categorical context:
\begin{prp} An exact functor $\psi\colon\fromto{\mathscr{B}}{\mathscr{A}}$ between exact $\infty$-categories induces an equivalence of $K$-theory spectra if, for every object $X\in\mathscr{B}$, the simplicial set
\begin{equation*}
Q(\psi)_{/X}\coloneq Q(\mathscr{B})\times_{Q(\mathscr{A})}Q(\mathscr{A})_{/X}
\end{equation*}
is weakly contractible.
\end{prp}

There is also the following variant of Quillen's Theorem B.
\begin{thm}[Theorem B for $\infty$-categories]\label{thm:B} Suppose $G\colon\fromto{C}{D}$ a functor between $\infty$-categories. If, for any morphism $f\colon\fromto{X}{Y}$ of $D$, the map
\begin{equation*}
f^{\star}\colon\fromto{G_{X/}}{G_{Y/}}
\end{equation*}
is a weak homotopy equivalence, then for any object $X\in D$, the square
\begin{equation*}
\begin{tikzpicture} 
\matrix(m)[matrix of math nodes, 
row sep=4ex, column sep=4ex, 
text height=1.5ex, text depth=0.25ex] 
{G_{X/}&C\\ 
D_{X/}&D\\}; 
\path[>=stealth,->,font=\scriptsize] 
(m-1-1) edge (m-1-2) 
edge (m-2-1) 
(m-1-2) edge (m-2-2) 
(m-2-1) edge (m-2-2); 
\end{tikzpicture}
\end{equation*}
is a homotopy pullback (for the Quillen model structure), and of course $D_{X/}\simeq\ast$.
\begin{proof} Consider the diagram
\begin{equation*}
\begin{tikzpicture} 
\matrix(m)[matrix of math nodes, 
row sep=4ex, column sep=4ex, 
text height=1.5ex, text depth=0.25ex] 
{G_{X/}&\widetilde{\mathscr{O}}(D)\times_{D}C&C\\ 
D_{X/}&\widetilde{\mathscr{O}}(D)&D\\
\{X\}&D^{\op}&\\}; 
\path[>=stealth,->,font=\scriptsize] 
(m-1-1) edge (m-1-2) 
edge (m-2-1) 
(m-1-2) edge (m-1-3)
edge (m-2-2)
(m-1-3) edge (m-2-3) 
(m-2-1) edge (m-2-2)
edge (m-3-1)
(m-2-2) edge (m-2-3)
edge (m-3-2)
(m-3-1) edge (m-3-2); 
\end{tikzpicture}
\end{equation*}
in which every square is a pullback. The upper right and lower left squares are homotopy pullbacks because opposite maps in these squares are weak homotopy equivalences. It therefore remains to show that the large left-hand rectangle is a homotopy pullback. For this, it is enough to show that $\fromto{\widetilde{\mathscr{O}}(D)\times_{D}C}{D^{\op}}$ satisfies the conditions Waldhausen's Theorem B for simplicial sets \cite[Lm. 1.4.B]{MR86m:18011} or, equivalently, is a sharp map in the sense of Hopkins and Rezk \cite[\S2]{math.AT/9811038}.

To prove this, suppose $\sigma\colon\fromto{\Delta^n}{D^{\op}}$ an $n$-simplex. Denote by $\sigma(0)$ its restriction to the $0$-simplex $\Delta^{\{0\}}\subset\Delta^n$. The pullback
\begin{equation*}
G_{\sigma/}\simeq(\widetilde{\mathscr{O}}(D)\times_{D}C)\times_{D^{\op}}\Delta^n
\end{equation*}
is naturally equivalent as an $\infty$-category to $G_{\sigma(0)/}$, and for any map $\eta\colon\fromto{\Delta^m}{\Delta^n}$, the pullback
\begin{equation*}
G_{\sigma\circ\eta/}\simeq(\widetilde{\mathscr{O}}(D)\times_{D}C)\times_{D^{\op}}\Delta^m
\end{equation*}
is naturally equivalent to $G_{\sigma(\eta(0))/}$. Our hypothesis is precisely that the morphism $\fromto{\sigma(\eta(0))}{\sigma(0)}$ induces a weak homotopy equivalence $\fromto{G_{\sigma(0)/}}{G_{\sigma(\eta(0))/}}$, whence $\eta$ induces an equivalence
\begin{equation*}
\equivto{(\widetilde{\mathscr{O}}(D)\times_{D}C)\times_{D^{\op}}\Delta^n\simeq G_{\sigma/}}{G_{\sigma\circ\eta/}\simeq(\widetilde{\mathscr{O}}(D)\times_{D}C)\times_{D^{\op}}\Delta^m},
\end{equation*}
as desired.
\end{proof}
\end{thm}

The various conditions required in Quillen's D\'evissage Theorem \cite[Th. 4]{MR0338129} are quite stringent, even for abelian For stable $\infty$-categories, they are simply unreasonable. For example, the inclusion of the stable $\infty$-category of bounded complexes of finite-dimensional $\FF_p$-vector spaces into the stable $\infty$-category of $p$-torsion bounded complexes of finitely generated abelian groups is not full, and there is no meaningful sense in which it is closed under the formation of subobjects.

Nevertheless, one can make use of the sort of filtrations Quillen employed in his D\'evissage Theorem. To this end, let us suppose that $\psi\colon\fromto{\mathscr{B}}{\mathscr{A}}$ is an exact functor between \emph{stable} $\infty$-categories. Our aim is to reduce the study of the map $\fromto{Q(\mathscr{B})}{Q(\mathscr{A})}$ induced by $\psi$ to the study of the weak homotopy type of a single $\infty$-category $Z(\psi)$. The key lemma that makes this reduction possible is the following, which is analogous to the proof of Quillen's D\'evissage Theorem \cite[Th. 4]{MR0338129}. 
\begin{lem}\label{lem:keydevlemma} For any object $U\in\mathscr{B}$ and any pushout square
\begin{equation*}
\begin{tikzpicture} 
\matrix(m)[matrix of math nodes, 
row sep=4ex, column sep=4ex, 
text height=1.5ex, text depth=0.25ex] 
{Y&X\\ 
0&\psi U\\}; 
\path[>=stealth,->,font=\scriptsize] 
(m-1-1) edge node[above]{$\phi$} (m-1-2) 
edge (m-2-1) 
(m-1-2) edge (m-2-2) 
(m-2-1) edge (m-2-2); 
\end{tikzpicture}
\end{equation*}
of $\mathscr{A}$, the functor $\phi_!\colon\fromto{Q(\psi)_{/Y}}{Q(\psi)_{/X}}$ induced by $\phi$ is a weak homotopy equivalence.
\begin{proof} We define a functor $\phi^{\star}\colon\fromto{Q(\psi)_{/X}}{Q(\psi)_{/Y}}$, functors
\begin{equation*}
\lambda_{\phi}\colon\fromto{Q(\psi)_{/X}}{Q(\psi)_{/X}}\textrm{\quad and\quad}\mu_{\phi}\colon\fromto{Q(\psi)_{/Y}}{Q(\psi)_{/Y}},
\end{equation*}
and natural transformations
\begin{eqnarray}
&\phi_{!}\circ\phi^{\star}\ \tikz[baseline]\draw[>=stealth,->,font=\scriptsize](0,0.5ex)--node[above]{$\alpha$}(0.5,0.5ex);\ \lambda_{\phi}\ \tikz[baseline]\draw[>=stealth,<-,font=\scriptsize](0,0.5ex)--node[above]{$\beta$}(0.5,0.5ex);\ \id_{Q(\psi)_{/X}}&\nonumber\\
&\phi^{\star}\circ\phi_!\ \tikz[baseline]\draw[>=stealth,->,font=\scriptsize](0,0.5ex)--node[above]{$\gamma$}(0.5,0.5ex);\ \mu_{\phi}\ \tikz[baseline]\draw[>=stealth,<-,font=\scriptsize](0,0.5ex)--node[above]{$\delta$}(0.5,0.5ex);\ \id_{Q(\psi)_{/Y}}.&\nonumber
\end{eqnarray}
These will exhibit $\phi^{\star}$ as a homotopy inverse of $\phi_{!}$.

For any object $T$ of $\mathscr{A}$, an $n$-simplex $(W,Z,g)$ of $Q(\psi)_{/T}$ may be said to consist of:
\begin{itemize}
\item[---] a diagram $W$ in $\mathscr{B}$ of the form
\begin{equation*}
\begin{tikzpicture} 
\matrix(m)[matrix of math nodes, 
row sep=4ex, column sep=4ex, 
text height=1.5ex, text depth=0.25ex] 
{W_{0n}&W_{1n}&\cdots&W_{(n-1)n}&W_{nn}\\ 
W_{0(n-1)}&W_{1(n-1)}&\cdots&W_{(n-1)(n-1)}&\\
\vdots&\vdots&\revddots&&\\
W_{01}&W_{11}&&&\\
W_{00}&&&&\\}; 
\path[>=stealth,->,font=\scriptsize] 
(m-1-1) edge (m-1-2) 
edge (m-2-1) 
(m-1-2) edge (m-1-3)
edge (m-2-2)
(m-1-3) edge (m-1-4)
(m-1-4) edge (m-1-5)
edge (m-2-4) 
(m-2-1) edge (m-2-2)
edge (m-3-1)
(m-2-2) edge (m-2-3)
edge (m-3-2)
(m-2-3) edge (m-2-4)
(m-3-1) edge (m-4-1)
(m-3-2) edge (m-4-2)
(m-4-1) edge (m-4-2)
edge (m-5-1); 
\end{tikzpicture}
\end{equation*}
in which every square is a pullback square,
\item[---] a sequence $Z$ of edges in $\mathscr{A}_{/T}$ of the form
\begin{equation*}
Z_0\to Z_1\to\cdots\to Z_n\to T,
\end{equation*}
and
\item[---] a diagram $g$ in $\mathscr{A}$ of the form
\begin{equation*}
\begin{tikzpicture} 
\matrix(m)[matrix of math nodes, 
row sep=4ex, column sep=4ex, 
text height=1.5ex, text depth=0.25ex] 
{Z_0&Z_1&\cdots&Z_n\\ 
\psi W_{0n}&\psi W_{1n}&\cdots&\psi W_{nn}\\}; 
\path[>=stealth,->,font=\scriptsize] 
(m-1-1) edge (m-1-2) 
edge node[left]{$g_0$} (m-2-1) 
(m-1-2) edge (m-1-3)
edge node[left]{$g_1$} (m-2-2)
(m-1-3) edge (m-1-4)
(m-1-4) edge node[right]{$g_n$} (m-2-4)
(m-2-1) edge (m-2-2)
(m-2-2) edge (m-2-3)
(m-2-3) edge (m-2-4); 
\end{tikzpicture}
\end{equation*}
in which every square is a pullback.
\end{itemize}
In this notation, the functor $\phi_!$ carries a simplex $(W,Z,g)\in Q(\psi)_{/Y}$ to $(W,Z,g)\in Q(\psi)_{/X}$, where by an abuse of notation, we write $Z$ also for the image of $Z$ in $\mathscr{A}_{/X}$.

Now we may define the remaining functors as follows.
\begin{itemize}
\item[---] Let $\phi^{\star}$ be the functor $\fromto{Q(\psi)_{/X}}{Q(\psi)_{/Y}}$ that carries a simplex $(W,Z,g)\in Q(\psi)_{/X}$ to
\begin{equation*}
(W,Z\times_XY,g\circ\pr_Z)\in Q(\psi)_{/Y},
\end{equation*}
where $(Z\times_XY)_i\coloneq Z_i\times_XY$, and $(\pr_Z)_i$ is the projection $\fromto{Z_i\times_XY}{Z_i}$.
\item[---] Let $\lambda_{\phi}$ be the functor $\fromto{Q(\psi)_{/X}}{Q(\psi)_{/X}}$ that carries a simplex $(W,Z,g)\in Q(\psi)_{/X}$ to
\begin{equation*}
(W\oplus U,Z,j_U\circ g)\in Q(\psi)_{/X},
\end{equation*}
where $(W\oplus U)_i\coloneq W_i\oplus U$, and $(j_U)_i$ is the inclusion $\into{\psi W_i}{\psi(W_i\oplus U)}$.
\item[---] Let $\mu_{\phi}$ be the functor $\fromto{Q(\psi)_{/Y}}{Q(\psi)_{/Y}}$ that carries a simplex $(W,Z,g)\in Q(\psi)_{/Y}$ to
\begin{equation*}
(W\oplus\Omega U,Z\oplus\psi\Omega U,g\oplus\id_{\psi\Omega U})\in Q(\psi)_{/Y}.
\end{equation*}
\end{itemize}

The composite $\phi_!\circ\phi^{\star}\colon\fromto{Q(\psi)_{/X}}{Q(\psi)_{/X}}$ clearly carries $(W,Z,g)\in Q(\psi)_{/X}$ to
\begin{equation*}
(W,Z\times_XY,g\circ\pr_Z)\in Q(\psi)_{/X}.
\end{equation*}
On the other hand, the canonical equivalence $Y\times_XY\simeq Y\oplus\psi\Omega U$ permits us to express the composite $\phi^{\star}\circ\phi_!\colon\fromto{Q(\psi)_{/Y}}{Q(\psi)_{/Y}}$ as the functor that carries $(W,Z,g)\in Q(\psi)_{/Y}$ to
\begin{equation*}
(W,Z\oplus\psi\Omega U,g\circ\pr_Z)\in Q(\psi)_{/Y}. 
\end{equation*}

Now we are prepared to define the natural transformations $\alpha,\beta,\gamma,\delta$.
\begin{itemize}
\item[---] The component at $(W,Z,g)$ of the natural transformation $\alpha$ is the diagram
\begin{equation*}
\begin{tikzpicture} 
\matrix(m)[matrix of math nodes, 
row sep=4ex, column sep=4ex, 
text height=1.5ex, text depth=0.25ex] 
{Z\times_XY&Z&X\\ 
\psi W&\psi(W\oplus U)&\\
\psi W;&&\\}; 
\path[>=stealth,->,font=\scriptsize] 
(m-1-1) edge node[above]{$\pr_Z$} (m-1-2) 
edge node[left]{$g\circ\pr_Z$} (m-2-1) 
(m-1-2) edge (m-1-3)
edge node[right]{$j_U\circ g$} (m-2-2) 
(m-2-1) edge node[below]{$j_U$} (m-2-2)
edge[-,double distance=1.5pt] (m-3-1); 
\end{tikzpicture}
\end{equation*}
the square is a pullback since 
\begin{equation*}
\begin{tikzpicture} 
\matrix(m)[matrix of math nodes, 
row sep=4ex, column sep=4ex, 
text height=1.5ex, text depth=0.25ex] 
{Z\times_XY&Z\\ 
Z&Z\oplus\psi U\\}; 
\path[>=stealth,->,font=\scriptsize] 
(m-1-1) edge (m-1-2) 
edge (m-2-1) 
(m-1-2) edge (m-2-2) 
(m-2-1) edge (m-2-2); 
\end{tikzpicture}
\end{equation*}
is so.
\item[---] The component at $(W,Z,g)$ of the natural transformation $\beta$ is the diagram
\begin{equation*}
\begin{tikzpicture} 
\matrix(m)[matrix of math nodes, 
row sep=4ex, column sep=4ex, 
text height=1.5ex, text depth=0.25ex] 
{Z&Z&X\\ 
\psi(W\oplus U)&\psi(W\oplus U)&\\
\psi U.&&\\}; 
\path[>=stealth,->,font=\scriptsize] 
(m-1-1) edge[-,double distance=1.5pt] (m-1-2) 
edge node[left]{$j_U\circ g$} (m-2-1) 
(m-1-2) edge (m-1-3)
edge node[right]{$j_U\circ g$} (m-2-2) 
(m-2-1) edge[-,double distance=1.5pt] (m-2-2)
edge (m-3-1); 
\end{tikzpicture}
\end{equation*}
\item[---] The component at $(W,Z,g)$ of the natural transformation $\gamma$ is the diagram
\begin{equation*}
\begin{tikzpicture} 
\matrix(m)[matrix of math nodes, 
row sep=4ex, column sep=4ex, 
text height=1.5ex, text depth=0.25ex] 
{Z\oplus\psi\Omega U&Z\oplus\psi\Omega U&Y\\ 
\psi(W\oplus\Omega U)&\psi(W\oplus\Omega U)&\\
\psi U.&&\\}; 
\path[>=stealth,->,font=\scriptsize] 
(m-1-1) edge[-,double distance=1.5pt] (m-1-2) 
edge node[left]{$g\oplus\id_{\psi\Omega U}$} (m-2-1) 
(m-1-2) edge (m-1-3)
edge node[right]{$g\oplus\id_{\psi\Omega U}$} (m-2-2) 
(m-2-1) edge[-,double distance=1.5pt] (m-2-2)
edge (m-3-1); 
\end{tikzpicture}
\end{equation*}
\item[---] Finally, the component at $(W,Z,g)$ of the natural transformation $\delta$ is the diagram
\begin{equation*}
\begin{tikzpicture} 
\matrix(m)[matrix of math nodes, 
row sep=4ex, column sep=4ex, 
text height=1.5ex, text depth=0.25ex] 
{Z&Z\oplus\Omega U&X\\ 
\psi W&\psi(W\oplus\Omega U)&\\
\psi W.&&\\}; 
\path[>=stealth,->,font=\scriptsize] 
(m-1-1) edge node[above]{$j_Z$} (m-1-2) 
edge node[left]{$g$} (m-2-1) 
(m-1-2) edge (m-1-3)
edge node[right]{$g\oplus\id_{\psi\Omega U}$} (m-2-2) 
(m-2-1) edge node[below]{$j_U$} (m-2-2)
edge[-,double distance=1.5pt] (m-3-1); 
\end{tikzpicture}
\end{equation*}
\end{itemize}
This completes the proof.
\end{proof}
\end{lem}

\begin{nul} It's worthwhile to note now some distinctions between our argument and Quillen's. Our argument is more complicated in some ways and easier in others. (We suspect that one should be able to prove a result for more general exact $\infty$-categories that contains all the possible complications and specializes both to the lemma above and to Quillen's result, but the additional complications are unnecessary for our work here.) The functor $\mu_{\phi}$ above, for instance, does not make an appearance in Quillen's argument. In effect, this is because on an ordinary abelian category, the loop space functor $\Omega$ coincides the constant functor at $0$. Furthermore, Quillen has to make use of closure properties of the full subcategory in order to check that suitable kernels and cokernels exist. In the stable setting, these issues vanish. Finally, in Quillen's case, $Q(\psi)_{/0}$ is visibly contractible, and the lemma above along with the existence of suitable filtrations imply the D\'evissage Theorem. In our situation, it is not always the case that $Q(\psi)_{/0}$ weakly contractible; this condition has to be verified separately. But we can find analogues of the filtrations sought by Quillen.
\end{nul}

\begin{dfn}\label{dfn:devissable} A filtration
\begin{equation*}
0=X_0\to X_1\to\cdots\to X
\end{equation*}
of an object $X\in\mathscr{A}$ will be said to be \textbf{\emph{$\psi$-admissible}} if the following conditions are satisfied.
\begin{enumerate}[(\ref{dfn:devissable}.1)]
\item The diagram above exhibits $X$ as the colimit $\colim_i X_i$.
\item For any integer $i\geq 1$, the cofiber $C_i\coloneq X_i/X_{i-1}$ is equivalent to $\psi U_i$ for some object $U_i\in\mathscr{B}$.
\item For any corepresentable functor $F\colon\fromto{\mathscr{A}}{\Kan}$, the diagram
\begin{equation*}
0=FX_0\to FX_1\to\cdots\to FX
\end{equation*}
exhibits $FX$ as the colimit $\colim_i FX_i$.
\end{enumerate}

We will say that $\psi$ is a \textbf{\emph{nilimmersion}} if every object of $\mathscr{A}$ admits a $\psi$-admissible filtration. 
\end{dfn}

\begin{exm}\label{exm:Lambdax} Suppose $\Lambda$ an $E_1$ ring, and suppose $x\in\pi_{\ast}\Lambda$ a homogeneous element of degree $d$. Suppose $f\colon\fromto{\Lambda}{\Lambda'}$ a morphism of $E_1$ rings such that one has a cofiber sequence
\begin{equation*}
\Lambda[d]\ \tikz[baseline]\draw[>=stealth,->,font=\scriptsize](0,0.5ex)--node[above]{$x$}(0.5,0.5ex);\ \Lambda\ \tikz[baseline]\draw[>=stealth,->,font=\scriptsize](0,0.5ex)--node[above]{$f$}(0.5,0.5ex);\ \Lambda'
\end{equation*}
of left $\Lambda$-modules. Assume that the induced functor $\fromto{\Mod^{\ell}_{\Lambda'}}{\Mod^{\ell}_{\Lambda}}$ preserves compact objects, so that it restricts to a functor $\fromto{\Perf^{\ell}_{\Lambda'}}{\Perf^{\ell}_{\Lambda}}$. If the multiplicative system $S\subset\pi_{\ast}\Lambda$ generated by $x$ satisfies the left Ore condition, then by \cite[Lm. 8.2.4.13]{HA}, the natural functor $\fromto{\mathbf{Perf}^{\ell}_{\Lambda'}}{\mathbf{Nil}_{(\Lambda,S)}^{\ell,\omega}}$ is a nilimmersion \cite[Pr. 11.15]{K1}.
\end{exm}

\begin{exm} As a subexample, when $p$ is prime and $\Lambda=\mathrm{BP}\langle n\rangle$, we have a nilimmersion
\begin{equation*}
\fromto{\mathbf{Perf}^{\ell}_{\mathrm{BP}\langle n-1\rangle}}{\mathbf{Nil}_{(\mathrm{BP}\langle n\rangle,S)}^{\ell,\omega}}
\end{equation*}
where $S\subset\pi_{\ast}\Lambda$ is the multiplicative system generated by $v_n$. As suggested by \cite[Ex. 11.16]{K1}, this particular nilimmersion is of particular import for a well-known conjecture of Ausoni--Rognes \cite[(0.2)]{MR1947457}.
\end{exm}

\begin{thm}[``Proto-d\'evissage'']\label{thm:protodev} Suppose $\psi$ a nilimmersion. Then the square
\begin{equation*}
\begin{tikzpicture} 
\matrix(m)[matrix of math nodes, 
row sep=4ex, column sep=4ex, 
text height=1.5ex, text depth=0.25ex] 
{Q(\psi)_{/0}&Q(\mathscr{B})\\ 
Q(\mathscr{A})_{/0}&Q(\mathscr{A})\\}; 
\path[>=stealth,->,font=\scriptsize] 
(m-1-1) edge (m-1-2) 
edge (m-2-1) 
(m-1-2) edge (m-2-2) 
(m-2-1) edge (m-2-2); 
\end{tikzpicture}
\end{equation*}
is a homotopy pullback (for the Quillen model structure), and of course $Q(\mathscr{A})_{/0}\simeq\ast$. In particular, the $\infty$-category $Q(\psi)_{/0}$ is weakly contractible just in case the induced map $\fromto{K(\mathscr{B})}{K(\mathscr{A})}$ is a weak equivalence.
\begin{proof} We employ our variant of Theorem B (Th. \ref{thm:B}). The conditions of this theorem will be satisfied once we check that for any object $X\in\mathscr{A}$, the $\infty$-category $Q(\psi)_{/X}$ is weakly equivalent to $Q(\psi)_{/0}$. For this, we apply Lm. \ref{lem:keydevlemma} to a $\psi$-admissible filtration
\begin{equation*}
0=X_0\to X_1\to\cdots\to X
\end{equation*}
to obtain a sequence of weak homotopy equivalences of $\infty$-categories
\begin{equation*}
Z(\psi)\simeq Q(\psi)_{/0}=Q(\psi)_{/X_0}\ \tikz[baseline]\draw[>=stealth,->,font=\scriptsize,inner sep=0.5pt](0,0.5ex)--node[above]{$\sim$}(0.5,0.5ex);\ Q(\psi)_{/X_1}\ \tikz[baseline]\draw[>=stealth,->,font=\scriptsize,inner sep=0.5pt](0,0.5ex)--node[above]{$\sim$}(0.5,0.5ex);\ \cdots\ \tikz[baseline]\draw[>=stealth,->,font=\scriptsize,inner sep=0.5pt](0,0.5ex)--(0.5,0.5ex);\ Q(\psi)_{/X}.
\end{equation*}
The result now follows from the claim that the diagram above exhibits $Q(\psi)_{/X}$ as the homotopy colimit (even in the Joyal model structure) of the $\infty$-categories $Q(\psi)_{/X_i}$.

To prove this claim, it is enough to show the following: (1) that the set $\pi_0\iota Q(\psi)_{/X}$ is exhibited as the colimit of the sets $\pi_0\iota Q(\psi)_{/X_i}$, and (2) that for any objects $A,B\in Q(\psi)_{/X_i}$, the natural map
\begin{equation*}
\fromto{\colim_{j\geq i}\Map_{Q(\psi)_{/X_j}}(A,B)}{\Map_{Q(\psi)_{/X}}(A,B)}
\end{equation*}
is an equivalence. The first claim follows directly, since any morphism $\fromto{Y}{X}$ factors through a morphism $\fromto{Y}{X_i}$ for some integer $i\geq 0$ by (\ref{dfn:devissable}.3). For the second, note that for any object $Y\in\mathscr{A}$, one may identify a mapping space $\Map_{Q(\psi)_{/Y}}(A,B)$ as the homotopy fiber of the map 
\begin{equation*}
\fromto{\Map_{Q(\psi)_{/0}}(A,B)}{\Map_{\mathscr{A}}(s(A),s(Y))}
\end{equation*}
induced the source functor $s\colon\fromto{Q(\psi)_{/0}}{\mathscr{A}}$ and the given map $\fromto{s(B)}{s(Y)}$ over the point corresponding to the given map $\fromto{s(A)}{s(Y)}$. So the natural map above can be rewritten as the natural map
\begin{equation*}
\begin{tikzpicture} 
\matrix(m)[matrix of math nodes, 
row sep=4ex, column sep=4ex, 
text height=1.5ex, text depth=0.25ex] 
{\colim_{j\geq i}\left(\Map_{Q(\psi)_{/0}}(A,B)\times_{\Map_{\mathscr{A}}(s(A),s(X_j))}\{\phi_j\}\right)\\ 
\Map_{Q(\psi)_{/0}}(A,B)\times_{\Map_{\mathscr{A}}(s(A),s(X))}\{\phi\};\\}; 
\path[>=stealth,->,font=\scriptsize] 
(m-1-1) edge (m-2-1); 
\end{tikzpicture}
\end{equation*}
where $\phi_j\colon\fromto{s(A)}{s(X_j)}$ and $\phi\colon\fromto{s(A)}{s(X)}$ are the given maps. Now since filtered homotopy colimits commute with homotopy pullbacks, the proof is complete thanks to (\ref{dfn:devissable}.3).
\end{proof}
\end{thm}

Let us study the $\infty$-category $Q(\psi)_{/0}$ more explicitly.

\begin{cnstr} The $\infty$-category $\mathscr{O}(\mathscr{A})\coloneq\Fun(\Delta^1,\mathscr{A})$ is also stable, and the cartesian and cocartesian fibration $t\colon\fromto{\mathscr{O}(\mathscr{A})}{\mathscr{A}}$ induces a cartesian fibration
\begin{equation*}
p\colon\fromto{Q(\mathscr{O}(\mathscr{A}))}{Q(\mathscr{A})}.
\end{equation*}
One observes that $p$ is classified by a functor $\fromto{Q(\mathscr{A})^{\op}}{\Cat_{\infty}}$ that carries an object $X\in\mathscr{A}$ to the $\infty$-category $\mathscr{A}_{/X}$. We may now extract the maximal right fibration
\begin{equation*}
\iota_{Q(\mathscr{A})}p\colon\fromto{\iota_{Q(\mathscr{A})}Q(\mathscr{O}(\mathscr{A}))}{Q(\mathscr{A})}
\end{equation*}
contained in $p$. Write
\begin{equation*}
E(\mathscr{A})\coloneq\iota_{Q(\mathscr{A})}Q(\mathscr{O}(\mathscr{A})).
\end{equation*}
In light of \cite[Cor. 3.3.4.6]{HTT}, the simplicial set $E(\mathscr{A})$ is the (homotopy) colimit of the functor $\fromto{Q(\mathscr{A})^{\op}}{\Kan}$ that classifies $\iota_{Q(\mathscr{A})}p$. This functor carries an object $X\in\mathscr{A}$ to the Kan complex $\iota(\mathscr{A}_{/X})$, and it carries a morphism
\begin{equation*}
\begin{tikzpicture} 
\matrix(m)[matrix of math nodes, 
row sep=3ex, column sep=3ex, 
text height=1.5ex, text depth=0.25ex] 
{&Z&\\ 
Y&&X\\}; 
\path[>=stealth,->,font=\scriptsize] 
(m-1-2) edge (m-2-1) 
edge (m-2-3); 
\end{tikzpicture}
\end{equation*}
to the map $\fromto{\iota(\mathscr{A}_{/X})}{\iota(\mathscr{A}_{/Y})}$ given by $\goesto{T}{T\times_XZ}$.

Using the uniqueness of limits and colimits in $\infty$-categories \cite[Pr. 1.2.12.9]{HTT}, we obtain a trivial fibration $\equivto{Q(\mathscr{A})_{/0}}{E(\mathscr{A})}$; hence the $\infty$-category $Q(\psi)_{/0}$ is naturally equivalent to the $\infty$-category
\begin{equation*}
Z(\psi)\coloneq E(\mathscr{A})\times_{Q(\mathscr{A})}Q(\mathscr{B}).
\end{equation*}
These equivalences are compatible with the maps to $Q(\mathscr{B})$ and $Q(\mathscr{A})$, so the square
\begin{equation*}
\begin{tikzpicture} 
\matrix(m)[matrix of math nodes, 
row sep=4ex, column sep=4ex, 
text height=1.5ex, text depth=0.25ex] 
{Z(\psi)&Q(\mathscr{B})\\ 
E(\mathscr{A})&Q(\mathscr{A})\\}; 
\path[>=stealth,->,font=\scriptsize] 
(m-1-1) edge (m-1-2) 
edge (m-2-1) 
(m-1-2) edge (m-2-2) 
(m-2-1) edge (m-2-2); 
\end{tikzpicture}\end{equation*}
is a homotopy pullback (for the Quillen model structure), and of course $E(\mathscr{A})\simeq\ast$

Again employing \cite[Cor. 3.3.4.6]{HTT}, we find that the simplicial set $Z(\psi)$ is the (homotopy) colimit of the functor $\fromto{Q(\mathscr{B})^{\op}}{\Kan}$ that classifies the pulled back right fibration $\iota_{Q(\mathscr{A})}p\times_{Q(\mathscr{A})}Q(\mathscr{B})$. This functor carries an object $U\in\mathscr{B}$ to the Kan complex $\iota(\mathscr{A}_{/\psi U})$, and it carries a morphism of the form
\begin{equation*}
\begin{tikzpicture} 
\matrix(m)[matrix of math nodes, 
row sep=3ex, column sep=3ex, 
text height=1.5ex, text depth=0.25ex] 
{&W&\\ 
V&&U\\}; 
\path[>=stealth,->,font=\scriptsize] 
(m-1-2) edge (m-2-1) 
edge (m-2-3); 
\end{tikzpicture}
\end{equation*}
to the map $\fromto{\iota(\mathscr{A}_{/\psi U})}{\iota(\mathscr{A}_{/\psi V})}$ given by $\goesto{T}{T\times_{\psi U}\psi W}$.
\end{cnstr} 

\begin{dfn} We call the $\infty$-category $Z(\psi)$ constructed above the \textbf{\emph{relative $Q$ construction}} for the nilimmersion $\psi\colon\fromto{\mathscr{B}}{\mathscr{A}}$.
\end{dfn}

\begin{nul} To unpack this further, we may think of the objects of the $\infty$-category $Z(\psi)$ as pairs
\begin{equation*}
(U,g)=(U,g\colon\fromto{X}{\psi U}),
\end{equation*}
where $U\in\mathscr{B}$, and $g$ is a map of $\mathscr{A}$. A morphism
\begin{equation*}
\fromto{(V,h\colon\fromto{Y}{\psi V})}{(U,g\colon\fromto{X}{\psi U})}
\end{equation*}
of this $\infty$-category is then a pair of diagrams
\begin{equation*}
\left(
\begin{tikzpicture}[baseline]
\matrix(m)[matrix of math nodes, 
row sep=4ex, column sep=4ex, 
text height=1.5ex, text depth=0.25ex] 
{W&U\\ 
V&\\}; 
\path[>=stealth,->,font=\scriptsize] 
(m-1-1) edge (m-1-2) 
edge (m-2-1); 
\end{tikzpicture}
,\ 
\begin{tikzpicture}[baseline]
\matrix(m)[matrix of math nodes, 
row sep=4ex, column sep=4ex, 
text height=1.5ex, text depth=0.25ex] 
{Y&X\\ 
\psi W&\psi U\\
\psi V&\\}; 
\path[>=stealth,->,font=\scriptsize] 
(m-1-1) edge (m-1-2) 
edge (m-2-1) 
(m-1-2) edge (m-2-2) 
(m-2-1) edge (m-2-2)
edge (m-3-1); 
\end{tikzpicture}
\right),
\end{equation*}
the first from $\mathscr{B}$ and the second from $\mathscr{A}$, in which the square in the second diagram is a pullback.
\end{nul}

\begin{nul} Waldhausen introduced a relative $S_{\bullet}$ construction, which we described as a virtual Waldhausen $\infty$-category $\mathscr{K}(\psi)$ in \cite[Nt. 8.8]{K1}. This has the property that the sequence
\begin{equation*}
K(\mathscr{K}(\psi))\ \tikz[baseline]\draw[>=stealth,->](0,0.5ex)--(0.5,0.5ex);\ K(\mathscr{B})\ \tikz[baseline]\draw[>=stealth,->](0,0.5ex)--(0.5,0.5ex);\ K(\mathscr{A})
\end{equation*}
is a (homotopy) fiber sequence. In effect, $\mathscr{K}(\psi)$ is the geometric realization of the simplicial Waldhausen $\infty$-category $\KK_{\ast}(\psi)$ whose $m$-simplices consist of a totally filtered object
\begin{equation*}
0\ \tikz[baseline]\draw[>=stealth,>->](0,0.5ex)--(0.75,0.5ex);\ U_1\ \tikz[baseline]\draw[>=stealth,>->](0,0.5ex)--(0.75,0.5ex);\ U_2\ \tikz[baseline]\draw[>=stealth,>->](0,0.5ex)--(0.75,0.5ex);\ \dots\ \tikz[baseline]\draw[>=stealth,>->](0,0.5ex)--(0.75,0.5ex);\ U_m
\end{equation*}
of $\mathscr{B}$, a filtered object
\begin{equation*}
X_0\ \tikz[baseline]\draw[>=stealth,>->](0,0.5ex)--(0.75,0.5ex);\ X_1\ \tikz[baseline]\draw[>=stealth,>->](0,0.5ex)--(0.75,0.5ex);\ X_2\ \tikz[baseline]\draw[>=stealth,>->](0,0.5ex)--(0.75,0.5ex);\ \dots\ \tikz[baseline]\draw[>=stealth,>->](0,0.5ex)--(0.75,0.5ex);\ X_m
\end{equation*}
of $\mathscr{A}$, and a diagram
\begin{equation*}
\begin{tikzpicture} 
\matrix(m)[matrix of math nodes, 
row sep=4ex, column sep=4ex, 
text height=1.5ex, text depth=0.25ex] 
{X_0&X_1&X_2&\dots&X_m\\ 
0&\psi(U_1)&\psi(U_2)&\dots&\psi(U_m)\\}; 
\path[>=stealth,->,font=\scriptsize] 
(m-1-1) edge[>->] (m-1-2) 
edge (m-2-1)
(m-1-2) edge[>->] (m-1-3)
edge (m-2-2)
(m-1-3) edge[>->] (m-1-4)
edge (m-2-3)
(m-1-4) edge[>->] (m-1-5)
(m-1-5) edge (m-2-5)
(m-2-1) edge[>->] (m-2-2) 
(m-2-2) edge[>->] (m-2-3)
(m-2-3) edge[>->] (m-2-4)
(m-2-4) edge[>->] (m-2-5); 
\end{tikzpicture}
\end{equation*}
of $\mathscr{A}$ in which every square is a pushout.

One may now point out that the relative $Q$ construction $Z(\psi)$ can be identified with the edgewise subdivision of (a version of) the simplicial space $\goesto{\mathbf{m}}{\iota\KK_m(\psi)_{0}}$. So, in these terms, the argument above essentially ensures that the natural map
\begin{equation*}
\fromto{I(\mathscr{K}(\psi))}{\Omega I(\mathscr{SK}(\psi))\simeq K(\mathscr{K}(\psi))}
\end{equation*}
is a weak homotopy equivalence, where $I$ is the left derived functor of $\iota$.
\end{nul}

\begin{exm} Keep the notations and conditions of Ex. \ref{exm:Lambdax}. In light of \cite[Pr. 11.15]{K1}, the resulting sequence
\begin{equation*}
\KK(\Lambda')\ \tikz[baseline]\draw[>=stealth,->](0,0.5ex)--(0.5,0.5ex);\ \KK(\Lambda)\ \tikz[baseline]\draw[>=stealth,->](0,0.5ex)--(0.5,0.5ex);\ \KK(\Lambda[x^{-1}])
\end{equation*}
is a fiber sequence if and only if the $\infty$-category
\begin{equation*}
Z(f_{\star})=E(\mathbf{Nil}_{(\Lambda,S)}^{\ell,\omega})\times_{Q(\mathbf{Nil}_{(\Lambda,S)}^{\ell,\omega})}Q(\mathbf{Perf}^{\ell}_{\Lambda'})
\end{equation*}
is weakly contractible.
\end{exm}

\begin{exm}\label{exm:ARconj} As a subexample, let us note that for a prime $p$, when $\Lambda=\mathrm{BP}\langle n\rangle$ is the truncated Brown--Peterson spectrum, we find, using \cite[Ex. 11.16]{K1}, that the conjecture of Ausoni--Rognes \cite[(0.2)]{MR1947457} that the sequence
\begin{equation*}
\KK(\mathrm{BP}\langle n-1\rangle)\ \tikz[baseline]\draw[>=stealth,->](0,0.5ex)--(0.5,0.5ex);\ \KK(\mathrm{BP}\langle n\rangle)\ \tikz[baseline]\draw[>=stealth,->](0,0.5ex)--(0.5,0.5ex);\ \KK(E(n)).
\end{equation*}
is a fiber sequence is equivalent to the weak contractibility of the $\infty$-category
\begin{equation*}
Z(f_{\star})=E(\mathbf{Nil}_{(\mathrm{BP}\langle n\rangle,v_n)}^{\ell,\omega})\times_{Q(\mathbf{Nil}_{(\mathrm{BP}\langle n\rangle,v_n)}^{\ell,\omega})}Q(\mathbf{Perf}^{\ell}_{\mathrm{BP}\langle n-1\rangle}).
\end{equation*}
For the sake of clarity, let us state explicitly that the objects of $Z(f_{\star})$ are pairs
\begin{equation*}
(U,g)=(U,g\colon\fromto{X}{f_{\star}U}),
\end{equation*}
in which $U$ is a perfect left $\mathrm{BP}\langle n-1\rangle$-module, and $g$ is a map of perfect $v_n$-nilpotent left $\mathrm{BP}\langle n\rangle$-modules. A morphism \begin{equation*}
\fromto{(V,h\colon\fromto{Y}{f_{\star}V})}{(U,g\colon\fromto{X}{f_{\star}U})}
\end{equation*}
of this $\infty$-category is then a pair of diagrams
\begin{equation*}
\left(
\begin{tikzpicture}[baseline]
\matrix(m)[matrix of math nodes, 
row sep=4ex, column sep=4ex, 
text height=1.5ex, text depth=0.25ex] 
{W&U\\ 
V&\\}; 
\path[>=stealth,->,font=\scriptsize] 
(m-1-1) edge (m-1-2) 
edge (m-2-1); 
\end{tikzpicture}
,\ 
\begin{tikzpicture}[baseline]
\matrix(m)[matrix of math nodes, 
row sep=4ex, column sep=4ex, 
text height=1.5ex, text depth=0.25ex] 
{Y&X\\ 
f_{\star}W&f_{\star}U\\
f_{\star}V&\\}; 
\path[>=stealth,->,font=\scriptsize] 
(m-1-1) edge (m-1-2) 
edge (m-2-1) 
(m-1-2) edge (m-2-2) 
(m-2-1) edge (m-2-2)
edge (m-3-1); 
\end{tikzpicture}
\right),
\end{equation*}
the first in perfect left $\mathrm{BP}\langle n-1\rangle$-modules and the second in perfect $v_n$-nilpotent left $\mathrm{BP}\langle n\rangle$-modules, in which the square in the second diagram is a pullback. The conjecture is simply that this $\infty$-category is weakly contractible. The relatively concrete nature of this formulation of the Ausoni--Rognes conjecture is tantalizing but, so far, we are embarrassed to report, frustrating.
\end{exm}


\bibliographystyle{amsplain}
\bibliography{kthy}

\end{document}